\documentclass[12pt]{article}

\usepackage[centertags]{amsmath}
\usepackage{color,amssymb}
\def\numberlikeadb{\global\def\theequation{\thesection.\arabic{equation}}}
\numberlikeadb

\newtheorem{Theorem}{Theorem}[section]
\newtheorem{thm}[Theorem]{Theorem}
\newtheorem{cor}[Theorem]{Corollary}
\newtheorem{Def}[Theorem]{Definition}
\newtheorem{lem}[Theorem]{Lemma}
 
\newtheorem{prop}[Theorem]{Proposition}
\newtheorem{rem}[Theorem]{Remark}

\def \A { {\mathscr{A}}}   \def\a {\alpha} \def \b {\beta} \def \bp {{\mbox{\boldmath$\pi$}}} \def \d { \delta} \def \D { \Delta}

\def\BG{{\mathscr B}(\G)}
\def\BH{{\mathscr B}({\mathscr H})}
\def \F {{\mathscr{F}}} \def \G {\Gamma}
\def \g {{\cal G}}
\def \H {{\mathscr{H}}}
\def \K {{\mathscr{K}}}
\def \L {{\mathscr{L}}}
  \def \l {\lambda}  \def \s {\sigma}
\def \t {\tilde}
\def \T {\Theta}
\def \U {\Upsilon}
 \def \z {\zeta} \def \Z {{\bf Z}}
 \def\PBDP {{\rm PBDP}}

 \def\( {\left(}  \def\) {\right)}
 \def \qed {\hfill $\Box$}

\def\eq{\begin{equation}}        \def\en{\end{equation}}
\def\eqa{\begin{eqnarray}}      \def\ena{\end{eqnarray}}
\def\eqs{\begin{eqnarray*}}     \def\ens{\end{eqnarray*}}

 \usepackage{mathrsfs,stmaryrd}
 \def\lb{{\lbag}}
 \def\rb{{\rbag}}

 \def\Ref#1{(\ref{#1})}
 \def\non{{\nonumber}}
 \def\tA{{\t\A}}
 \def\tG{{\t\G}}
 \def \tH {{\t\H}}

\def\m{{\mathscr{M}_\g\circ}}

\newcommand{\hXi}{{\hat\Xi}}
\newcommand{\XiyAy}{{\Xi_y\vert_{A_y}}}

\newcommand{\XiyAyc}{{\Xi_y\vert_{A_y^c}}}
\newcommand{\XiAyc}{{\Xi\vert_{A_y^c}}}
\newcommand{\XixAx}{{\Xi_x\vert_{A_x}}}
\newcommand{\XiAx}{{\Xi\vert_{A_x}}}
\newcommand{\XixAxc}{{\Xi_x\vert_{A_x^c}}}
\newcommand{\XiAxc}{{\Xi\vert_{A_x^c}}}

\newcommand{\XiyByc}{{\Xi_y\vert_{B_y^c}}}
\newcommand{\XiByc}{{\Xi\vert_{B_y^c}}}
\newcommand{\XixBx}{{\Xi_x\vert_{B_x}}}
\newcommand{\XiBx}{{\Xi\vert_{B_x}}}
\newcommand{\XixBxc}{{\Xi_x\vert_{B_x^c}}}
\newcommand{\XiBxc}{{\Xi\vert_{B_x^c}}}

\newcommand{\hXiByc}{{\hXi\vert_{B_y^c}}}

\newcommand{\hXiBxc}{{\hXi\vert_{B_x^c}}}

\newcommand{\XixyAxy}{{\Xi_{xy}\vert_{A_{xy}}}}
\newcommand{\XiAxy}{{\Xi\vert_{A_{xy}}}}

\newcommand{\XixyAxyc}{{\Xi_{xy}\vert_{A_{xy}^c}}}
\newcommand{\XiAxyc}{{\Xi\vert_{A_{xy}^c}}}
\newcommand{\XixyBxyc}{{\Xi_{xy}\vert_{B_{xy}^c}}}
\newcommand{\XiBxyc}{{\Xi|_{B_{xy}^c}}}

\newcommand{\hXiBxyc}{{\hXi\vert_{B_{xy}^c}}}

\def\lt{\l^{[2]}}

\def\btau {\mbox{\boldmath$\tau$}}

\def\nexto{\kern -0.54em}
\newcommand{\E}{{\rm {I\ \nexto E}}}
\def\prob{{\rm {I\ \nexto P}}}
\def\var{{\rm Var}}
\def\real{{\rm {I\ \nexto R}}}

\def\Ceka{\v Cekan\-avi\v cius}
\begin{document}

\title{Stein's method and locally dependent point process approximation}

\author{ Aihua Xia\footnote{Postal address: Department of Mathematics and Statistics,
the University of Melbourne,
VIC 3010, Australia. E-mail address: xia@ms.unimelb.edu.au}
 \ and \  Fuxi Zhang\footnote{Postal address: School of Mathematical Sciences, Peking University, Beijing 100871,
China. E-mail address: zhangfxi@math.pku.edu.cn}\\
The University of Melbourne and Peking University
}
\date{{{21 February, 2011}}}

\maketitle

\begin{abstract}
Random events in space and time often exhibit a locally dependent
structure. When the events are very rare and dependent structure is
not too complicated, various studies in the literature have shown
that Poisson and compound Poisson processes can provide adequate
approximations. However, the accuracy of approximations does not
improve or may even deteriorate when the mean number of events
increases. In this paper, we investigate an alternative family of
approximating point processes and establish Stein's method for their
approximations. We prove two theorems to accommodate respectively
the positively and negatively related dependent structures. Three
examples are given to illustrate that our approach can circumvent
the technical difficulties encountered in compound Poisson process
approximation [see Barbour \& M{\aa}nsson~(2002)] and our
approximation error bound decreases when the mean number of the
random events increases, in contrast to increasing bounds for
compound Poisson process approximation.



\vskip12pt \noindent\textit{Running title\/}: Locally dependent point process approximation
\end{abstract}

\maketitle

\section{Introduction}

\setcounter{equation}{0} Random events in space and time often
exhibit a locally dependent structure. When the events are very rare
and the dependent structure is not too complicated, a natural
approach is to declump the events into clusters then approximate the
positions of the clusters by a suitable Poisson process and the
sizes of the clusters by independent and identically distributed
random elements, as well documented in Aldous~(1989). Consequently,
compound Poisson and marked Poisson processes are often widely
accepted as the `best approximate models' for clustered rare events.

The first attempt to estimate the errors of Poisson process
approximation seems to go back to Brown~(1983) with errors measured
in the total variation distance, while the errors in the
L\'evy-Prohorov distance were not studied until Jacod \& Mano~(1988)
and Nikunen \& Valkeila~(1991) [see also Xia~(1993)]. All these
studies are based on the stochastic calculus approach with a
filtration, a compensator and coupling techniques as the tools to
quantify the distances. Barbour and Brown~(1992), clearly inspired
by the success of Stein's method in multivariate Poisson
approximation [Barbour~(1988)], laid down a general framework for
using Stein's method to estimate the Poisson process approximation
errors. Their framework can be well adjusted for errors expressed in
terms of Janossy densities, Palm distributions and compensators [see
Barbour, Brown \& Xia~(1998) and Xia~(2005)]. In terms of compound
Poisson process approximation, there seems no major advance until
Arratia, Goldstein \& Gordon~(1989) who replaced the original point
process with a new one carrying the information of locations and
cluster sizes separately so that the Stein-Chen method for Poisson
approximation can be employed to obtain useful error bounds. There
are enormous advantages for this approach if one can successfully
declump the point process, but the procedure of declumping is far
from obvious in applications. By contrast, Barbour \&
M{\aa}nsson~(2002) avoided declumping totally by setting a framework
of Stein's method so that the quality of approximation can be
studied directly, and the authors summarized that the direct
approach `has conceptual advantages, but entails technical
difficulties' in p.~1492. One of the main difficulties is that
Stein's factors, like their counterparts for compound Poisson random
variable approximation [see Barbour, Chen \& Loh~(1992), Barbour \&
Utev~(1998) and Barbour \& Utev~(1999)], are generally too crude to
use unless more conditions are imposed such as the compound Poisson
process is very close to a Poisson process. An immediate consequence
is that the error bounds obtained often deteriorate when the mean of
the point process increases, i.e., more information is available. On
the other hand, using the improved estimates for Stein's factors for
Poisson process approximation in Xia~(2005) [cf Brown, Weinberg \&
Xia~(2000)], Chen \& Xia~(2004) managed to produce error estimates
for Poisson process approximation to short range dependent rare
events and the estimates will remain small (but not improve either)
when the average number of events increases.

It is well-known that the central limit theorem often exhibits the
{\it large sample property}, i.e. the larger the sample size, the
better the approximation, as evidenced by the Berry--Esseen bound
[see Chen and Shao~(2004)]. If we are interested in the total counts
of rare and weakly dependent events, the Poisson law of small
numbers is the cornerstone of the area. However, the Poisson
approximation error does not enjoy the large sample property when
more rare events are counted [Barbour \& Hall~(1984)]. The
shortcoming is due to the fact that a Poisson distribution has only
one parameter to fiddle with while a normal distribution has two
parameters. When more parameters are introduced, this property can
be recovered [see Presman~(1983), Kruopis~(1986), \Ceka~(1997),
Barbour \& Xia~(1999), Brown \& Xia~(2001), R\"ollin~(2005)]. In
fact, Brown \& Xia~(2001) discovered a large family of distributions
that can achieve the same purpose.

The success of compound Poisson process approximation essentially
hinges on the fact that the events are very rare. It is tempting to
ask whether the approximation theory is still valid when the events
are less rare, more heavily dependent and the mean number of events
increases? One way to tackle this problem is to keep the
approximating process as a Poisson process but weaken the metric for
quantifying the difference between point processes [Schuhmacher \&
Xia~(2008)]. The weaker metric will naturally limit its
applicability. The second approach is to introduce more parameters
into the approximating point process models. To put the idea in
practice, Xia \& Zhang~(2008) introduced a family of point process
counterparts of approximating distributions suggested in Brown \&
Xia~(2001), and named them as the polynomial birth-death point
processes, or {\it PBDP} in short. In particular, Xia \&
Zhang~(2008) bounded the distance between the Bernoulli process with
a constant success probability and a suitable \PBDP\ in terms of the
Barbour-Brown distance (defined in section~\ref{steinmethod} below,
see also Barbour \& Brown~(1992)). The assumption of the constant
success probability plays the crucial role there because the
symmetric structure enables the authors to construct a suitable
coupling to directly compare the two distributions. The pilot study
shows that, for the Bernoulli process with the same success
probability, it is possible to recover the large sample property for
\PBDP\ approximation. The purpose of this paper is to
demonstrate that the
large sample property prevails among a large group of point
processes when these \PBDP\ are used as approximating models. To
this end, we set up the Stein equation of \PBDP\ approximation and
establish its Stein factors so that one can directly estimate the
difference between the distribution of a general point process and
that of a \PBDP.

Our paper is arranged as follows. In section~\ref{steinmethod}, we
briefly review the polynomial birth-death point processes introduced
in Xia \& Zhang~(2008),  lay down a foundation of Stein's method for
their approximation and conclude the section with estimates of
Stein's factors in terms of the Barbour-Brown metric. To make our
paper reader-friendly, we postpone the technical proofs of Stein's
factors to section~\ref{proofofsteinfactor}.
Section~\ref{locallydeptpointprocesses} is devoted to point
processes with locally dependent structures which are analogous with those
in Chen \& Shao~(2004). We state two theorems for error
estimates of \PBDP\ approximations, respectively for positively and
negatively related dependence. The proofs of these theorems are
rather complicated so we leave them to the last two sections
(sections~\ref{prooftheoremld1} and \ref{prooftheoremld2}) of the
paper. Examples are provided in section~\ref{applications} to
illustrate the key steps of applying the main theorems.

\section{Stein's method for polynomial birth-death point processes}
\label{steinmethod}
\setcounter{equation}{0}

The family of approximating distributions in Brown \& Xia~(2001) was
introduced through the invariant distributions of birth-death
processes. For ease of use, they focused on the birth and death
rates as the polynomial functions of the states of the process, and
consequently called the invariant distribution as {\it polynomial
birth-death distribution}. More precisely, let
 \begin{equation} \label{bdrate}
\a_k = a + b k, \ \forall \ k \ge 0;  \ \ \ \b_k = k + \b k (k-1), \
\forall \ k \ge 0 ,
 \end{equation}
where $a> 0$, $0 \le b < 1$, $\b \ge 0$. A birth-death process with
birth rates $\{\a_k\}$ and death rates $\{\b_k\}$ must be ergodic.
As in Brown \& Xia~(2001), we let $Z_n(\cdot):=\{Z_n(t) : t \ge 0
\}$ be such a process with initial value $n$ and use $\pi_{a,b;\b}$ or simply $\pi$ when there is no confusion to stand for the invariant
distribution.

Let $\G$ be a compact metric space with metric $d_0$ bounded by 1
and Borel $\sigma$-algebra $\BG$ generated by $d_0$. Set $U, U_1,
U_2, \cdots$ as independent and identically distributed $\G$-valued
random elements with distribution $\mu$. In this paper, the
expression $\sum_{i=1}^X \d_{U_i}$ always implies that the
nonnegative integer random variable $X$ is independent of $\{ U_i :
\ i \ge 1\}$. We call $\Z$ a {\it polynomial birth-death point
process} [see Xia \& Zhang~(2008)] if it can be expressed as
 $$
\Z = \sum_{i=1}^Z \d_{U_i}
 $$
for $Z \sim \pi_{a,b;\b}$, and denote $\L (\Z)$ by
$\bp_{a,b;\b;\mu}$ or simply $\bp$ when there is no confusion. We
now give a few examples to illustrate that the definition is a
natural extension of the polynomial birth-death distribution.

\noindent{\bf Example 1} Suppose $Z$ follows ${\rm Binomial}(n,p)$,
then $\Z$ reduces to a binomial
process.

\noindent{\bf Example 2} If $Z$ is a Poisson random variable with mean $a$, then $\Z$ becomes a Poisson process on $\G$ with mean measure
$a \mu$.

\noindent{\bf Example 3} When $Z$ has a negative binomial distribution, we call $\Z$ a {\it negative binomial
process}.

 \begin{rem}
{\rm There are two possible ways to define a negative binomial
process. The one we defined here does not have the property of
independent increments while if we define it as a compound Poisson
process with clusters following a logarithmic distribution, then it
does have the property of independent increments. Nevertheless, the
two distributions converge when the intensity of the Poisson
component becomes large [see Remark~\ref{compp3} below].}
 \end{rem}

Now we construct a Markov process with invariant distribution
$\bp = \bp_{a,b;\b;\mu}$. Allowing repeats of points, each
finite integer-valued measure on $\G$ can be written as $\xi =
\sum_{i=1}^n \d_{x_i}$. Since the points $x_1,\cdots, x_n$ are not
necessarily distinct, we introduce the notation $\lb x_1 ,\cdots,
x_n \rb$ to stand for the collection of the $n$ points. In this
paper, we do not distinguish $\sum_{i=1}^n \d_{x_i}$ with the
collection $\lb x_1 ,\cdots, x_n \rb$, or a configuration with $n$
particles respectively located at $x_1, \cdots, x_n$. For example,
when we say a site/point $x$ or a particle at $x$ in $\xi$, it means
that $\xi (\{x\}) \ge1$.

For each measure $\xi$ on $\G$, we denote its total mass by $|\xi|$.
Let $\H$ be the class of all possible finite integer-valued measures
(also known as the configurations of point processes) on $\G$ and
let $\BH$ be the smallest $\sigma$-algebra in $\H$ making the
mappings $\xi\mapsto \xi(C)$ measurable for all relatively compact
Borel sets $C\subset\Gamma$. For each suitable measurable function
$h$ on $\H$, we define
 \begin{eqnarray}
\A h (\xi)
 & := &
\big( a + b |\xi| \big)  \int_\G \big( h (\xi + \d_x) - h (\xi)
\big) \mu (dx) \non
 \\ & &
+ \big( 1 + \b( |\xi| - 1) \big) \int_\G \big( h(\xi - \d_x) -
h(\xi) \big) \xi (d x) \non
 \\ & = &
\big( a + b |\xi| \big) \( \E h (\xi + \d_U) - h (\xi)  \) \non\\
&&+  \big(
1 + \b( |\xi| - 1) \big) \( \E h (\xi - \d_{V(\xi)}) - h (\xi) \) ,
\label{generator}
 \end{eqnarray}
where, for $\xi = \sum_{i=1}^n \d_{x_i}$, $V(\xi)$ is a uniformly
distributed random element on the collection $\lb x_1 ,\cdots,x_n
\rb$. In other words, $V(\xi)$ is equally likely to be one of $x_1$,
$\dots$, $x_n$. A particle system $\Z_\xi (\cdot):=\{ \Z_\xi (t) : t
\ge 0 \}$ with the generator $\A$ evolves as follows:
 \begin{itemize}
 \item
with rate $a$ a new particle immigrates to $\G$ and settles at a site
according to $\mu$;
 \item
with rate $b$ an existing particle gives a birth, and the new born
particle is also located at a site chosen according to $\mu$;
 \item
with rate 1, an existing particle suicides;
 \item
with rate $\b$, an existing particle kills another existing
particle.
 \end{itemize}
We call such a Markov process as a {\it birth-death system}. It's
not difficult to check that the birth-death system has the unique
invariant distribution $\bp_{a,b;\b;\mu}$. Noting that
for any $\xi \in \H$, $\{ |\Z_\xi(t)| : t \ge 0 \}$ is a birth-death
process with rates (\ref{bdrate}), we have $\L(|\Z_\xi (\cdot)|) =
\L (Z_{|\xi|}(\cdot))$. Therefore, $\L\left(\Z_\xi (t) \right)=
\L\left(\sum_{i=1}^{Z_n (t)} \d_{U_i}\right)$ if $\L\xi
=\L\left(\sum_{i=1}^n \d_{U_i}\right)$. In particular, we have
$\L\left(\Z_\emptyset (t)\right)=\L\left(\sum_{i=1}^{Z_0 (t)}
\d_{U_i}\right)$.

Bearing in mind the Stein equation suggested by Barbour \& Brown~(1992), the natural choice of the Stein equation for the generator $\A$ is
\begin{equation} \label{Ah}
\A h(\xi) = f(\xi) -  \bp (f)
 \end{equation}
for suitable functions $f$ on $\H$, where $\bp(f):=  \int
f(\xi)\bp(d \xi) $. We now consider the question of the existence
of an $h$ that solves the equation \Ref{Ah}.

 \begin{prop} \label{solution}
For any bounded function $f$ on $\H$,
 $$
h_f (\xi) := - \int_0^\infty \big( \E f (\Z_\xi (t)) - \bp (f)
\big) d t
 $$
is well defined, and is a solution of (\ref{Ah}).
 \end{prop}

\noindent{\it Proof.} Let $\{U_i\}$ be independent $\mu$-distributed
random elements which are independent of $\{ \Z_\xi (t) : t \ge
0\}$. Pair $\{U_i,\ 1\le i\le |\xi|\}$ with the points in $\xi$,
define $\xi^\prime= \sum_{i=1}^{|\xi|} \d_{U_i}$, and construct $\{
\Z_{\xi^\prime} (t) : t \ge 0 \}$ from $\{ \Z_\xi (t) : t \ge 0\}$
by replacing the points in $\xi$ with the paired counterparts in
$\xi^\prime$. Let $\tilde{\tau}$ be the last death time of all the
points in $\xi$. We have
 $$
\int_0^\infty | \E f(\Z_\xi (t)) - \E f (\Z_{\xi^\prime} (t)) | d t
\le  \int_0^\infty \E \big( 2 \| f \|  1_{\tilde{\tau}>t} \big)
dt  = 2 \| f \| \E \tilde{\tau} < \infty,
 $$
since $\tilde{\tau}$ is stochastically smaller than the maximum of
$|\xi|$ independent and identically distributed $\exp(1)$ random variables.

Next, define $\bar{f} (n) = \E f (\sum_{i=1}^n \d_{U_i})$ for all $n\ge0$, then
$$\int_0^\infty \left|\E f (\Z_{\xi^\prime} (t)) - \bp (f) \right| d t
  \le
\int_0^\infty \left|\E \bar{f} (Z_{|\xi|} (t)) - \pi
(\bar{f})\right|d t<\infty
 $$
due to the positive recurrence of the Markov chain $\{Z_{|\xi|} (t),\ t\ge 0\}$. Hence,
\begin{eqnarray*}
&&\int_0^\infty | \E f(\Z_\xi (t)) -  \bp (f) | d
t\\
&& \le \int_0^\infty | \E f(\Z_\xi (t)) - \E f (\Z_{\xi^\prime} (t)) |
d t+\int_0^\infty \left|\E f (\Z_{\xi^\prime} (t)) - \bp (f) \right|
d t<\infty ,
 \end{eqnarray*}
which implies that $h_f$ is well-defined.

To establish \Ref{Ah}, let $\btau_\xi = \inf \{t : \Z_\xi (t) \neq \xi \}$, which has an exponential distribution with parameter $\a_{|\xi|} + \b_{|\xi|}$. Then
 \begin{eqnarray*}
h_f (\xi)
 & = &
- \int_0^\infty \big( \E f (\Z_\xi (t)) -  \bp(f) \big) d t
 \\ & = &
-  \big( f (\xi) -  \bp(f) \big) \E \btau_\xi - \E \int_{\btau_\xi}^\infty \big( \E f (\Z_\xi (t)) -  \bp(f) \big) d t
 \\ & = &
- \frac{ f (\xi) -  \bp(f) }{\a_{|\xi|} + \b_{|\xi|} } +  \E h (\Z_\xi (\btau_\xi))
  \\ & = &
- \frac{ f (\xi) -  \bp(f) } {\a_{|\xi|} + \b_{|\xi|} } +  \frac { \a_{|\xi|} \int_\G h ( \xi + \d_x) \mu(dx) + \big(1 + \b( |\xi| - 1) \big) \int_\G  h (\xi - \d_x) \xi(dx) } { \a_{|\xi|} + \b_{|\xi|} } ,
 \end{eqnarray*}
and \Ref{Ah} follows by rearranging the above equation. \qed

The metric used for quantifying the differences of two point
processes is defined as follows [see Barbour \& Brown~(1992)]. Let
$\K$ be the class of $d_0$-Lipschitz functions $u$ on $\G$ such that
$|u(x)-u(y)|\le d_0(x,y)$ for all $x,y\in\G$. For any two measures
$\rho_1$ and $\rho_2$ on $\G$, define
$$d_1(\rho_1,\rho_2)=\left\{\begin{array}{ll}
0,&\mbox{ if }|\rho_1|=|\rho_2|=0,\\
\frac1{|\rho_1|}\sup_{u\in\K}\left\vert\int_\G ud\rho_1-\int_\G u d\rho_2\right|,&\mbox{ if }|\rho_1|=|\rho_2|\ne 0,\\
1,&\mbox{ if }|\rho_1|\ne|\rho_2|.\end{array}\right.$$
For any configurations $\xi = \sum_{i=1}^n \d_{x_i}$ and $\eta = \sum_{i=1}^n
\d_{y_i} \in \H$ with $n\ge 1$, $d_1 (\xi, \eta)$ can be represented as
 $$
d_1 (\xi, \eta)=\min_\sigma \frac 1 n \sum_{i=1}^n d_0 (x_i,
y_{\sigma(i)}) ,
 $$
where the minimum is taken over all permutations $\s$ of $(1,\dots,n)$. The Barbour-Brown metric
$d_2$ between point process distributions is defined as
 $$
d_2 ({\bf P}, {\bf Q}) : = \sup_f |{\bf P} (f) - {\bf Q} (f)|
= \inf_{\xi \sim {\bf P}, \eta \sim {\bf Q}} \E d_1 (\xi, \eta) ,
 $$
where the supremum is taken over all functions in
 $$
\F: = \{ f: |f(\xi) - f(\eta)| \le d_1 (\xi, \eta), \ \forall \
\xi, \eta \in \H \} ,
 $$
and the last equation is due to the duality theorem [see
Rachev~(1991), p.~168]. The metric $d_2$ is a particular kind of the
well-known family of Wasserstein metrics. It is worthwhile to point
out that, since $d_1 \le 1$, all functions in $\F$ are bounded and
Proposition~\ref{solution} ensures the existence of solutions of
Stein's equation \Ref{Ah} for these functions. Historically, the
Wasserstein metrics were motivated by the classical
Monge-Transportation problem. In our context, we will handle the
`transportation problem' in two steps, i.e. to form `sandpiles'
by assembling local points to designated centers and then
transport the `sandpiles' of the point process being
approximated to the corresponding `sandpiles' of the PBDP.

The following Lemma is often useful for comparing two different approximating polynomial birth-death point processes.

\begin{lem}\label{lem2pbd1} We have
$$d_2\left(\bp_{a_1,b_1;\b_1;\mu_1},\bp_{a_2,b_2;\b_2;\mu_2}\right)\le d_{tv}(\pi_{a_1,b_1;\b_1},\pi_{a_2,b_2;\b_2})+d_1(\mu_1,\mu_2),$$
where for two probability measures $Q_1$ and $Q_2$ on
$\mathbb{Z}_+:=\{0,1,2,\dots\}$,
$$d_{tv}(Q_1,Q_2) :=\sup_{A\subset\mathbb{Z}_+}|Q_1(A)-Q_2(A)|.$$
 \end{lem}

\noindent{\it Proof.}  Using the Kantorovich-Rubinstein duality theorem [Rachev~(1991), Theorem~8.1.1,
p.~168], we can couple together $Z_1\sim\pi_{a_1,b_1;\b_1}$, $Z_2\sim\pi_{a_2,b_2;\b_2}$, and two sequences of $\G$-valued random elements $\tau_{1i}\sim\mu_1$ and $\tau_{2i}\sim\mu_2$, $i\ge 1$, such that
\eqs
&& d_{tv} \left(\pi_{a_1,b_1;\b_1},\pi_{a_2,b_2;\b_2}\right)=\prob(Z_1\ne Z_2),\\
&&\E d_0(\tau_{1i},\tau_{2i})=d_1\left(\mu_1,\mu_2\right) \mbox{ for all }i\ge1,
\ens
and $\{(\tau_{1i},\tau_{2i}),\ i\ge 1\}$ are independent and independent of $(Z_1,Z_2)$. Then
\eqs
&&d_2\left(\bp_{a_1,b_1;\b_1;\mu_1},\bp_{a_2,b_2;\b_2;\mu_2}\right)\le \E d_1\left(\sum_{i=1}^{Z_1} \d_{\tau_{1i}},\sum_{i=1}^{Z_2} \d_{\tau_{2i}}\right)\\
&\le&\prob(Z_1\ne Z_2)+ \E \left\{d_1\left.\left(\sum_{i=1}^{Z_1} \d_{\tau_{1i}},\sum_{i=1}^{Z_2} \d_{\tau_{2i}}\right)\right\vert Z_1=Z_2\right\}\prob(Z_1=Z_2)\\
&\le&d_{tv}\left(\pi_{a_1,b_1;\b_1},\pi_{a_2,b_2;\b_2}\right)+\E \left\{\left. {\frac 1 {Z_1} \sum_{i=1}^{Z_1} d_0\left(\tau_{1i},\tau_{2i}\right) }\right\vert Z_1=Z_2\right\}\prob(Z_1=Z_2)\\
&\le&d_{tv}(\pi_{a_1,b_1;\b_1},\pi_{a_2,b_2;\b_2})+d_1(\mu_1,\mu_2),
\ens completing the proof. \qed

In applications of Stein's equation, one will encounter the following quantities:
 \begin{equation} \label{Cn}
C_n : = \sup \{  | h_f (\xi + \d_x) - h_f (\xi + \d_y) |  : f \in
\F, \xi \in \H, |\xi| = n\} ,
 \end{equation}
with $C_{-1} := 0$,
$$
\D_2h (\xi;x,y) := h(\xi+\d_x+\d_y) - h(\xi +\d_x) - h(\xi+\d_y)
+h(\xi), \ \xi \in \H,\ x,y \in \G , $$
and
 $$
\D_2 h (\xi) := \sup \{ |\D_2 h(\xi; x, y)| :\ x,y \in \G \}.
 $$
The following estimates,  often known as Stein's factors, are
usually needed in applying Stein's method. If fact, the success of
Stein's method is centered around the quality of these estimates.

 \begin{thm} \label{mainresult}
(i) For $n \ge 0$,
 \begin{equation}
C_n \le \min \left\{ 1, \frac{1}{2(n+1)} + \frac 1 a , \frac{1}{ (a
\wedge b) (n+1)}  \right\} .\label{estimateCn}
\end{equation}

(ii) For any $f \in \F$, $\xi \in \H$,
\begin{equation}
\D_2 h_f (\xi) \le \frac 2 {|\xi|+1} + \frac 5 a.\label{mainresult2}
 \end{equation}
 \end{thm}

\begin{rem}
{\rm The estimates in Theorem~\ref{mainresult} are of the correct order. In fact, if we take $\b=b=0$, the \PBDP\ becomes a Poisson process and the estimates for the Poisson process are known to be of the correct order [see Xia~(2005)].}
 \end{rem}

\section{Locally dependent point processes}
\label{locallydeptpointprocesses}
\setcounter{equation}{0}

A {\it point process} $\Xi$ on $\G$ is defined as a measurable
mapping of some fixed probability space into $(\H,\BH)$ and
$\l(dx)=\E\Xi(dx)$ is said to be the {\it intensity} or {\it mean
measure} of $\Xi$ [Kallenberg~(1983), pp.~13-14]. A point process is
said to be {\it simple} if it has at most one point at each
location. For a point process $\Xi$ on $\G$ with finite mean measure
$\l$, the family of point processes $\{ \Xi_x : x \in \G \}$ are
said to be {\it reduced Palm processes} associated with $\Xi$ (at $x
\in\G $) if for any measurable function
$f:\G\times\H\rightarrow\real_+:=[0,\infty)$,
 \eq \label{eq:Palm}
\E\left(\int_{\G} f(x,\Xi-\d_x)\Xi(dx)\right) = \int_{\G} \E
f(x,\Xi_x) \l(dx) ,
 \en
[Kallenberg~(1983), Chapter 10]. Intuitively, the reduced Palm
distribution $\L\Xi_x$ is defined through the Radon-Nikodym
derivative as follows:
$$\prob(\Xi_x\in B)=\frac{\E[\Xi(dx)1_{\{\Xi-\d_x\in
B\}}]}{\E\Xi(dx)},\mbox{ for all }B\in\BH .
$$
When $\Xi$ is a simple point process, it can be interpreted as
the distribution of $\Xi$ save one point at $x$ conditional on there
is one point at $x$.

In this paper, we also need the {\it second order reduced Palm
processes} $\Xi_{xy}$ of the point process $\Xi$ at $x,\ y\in\G$
defined as the processes satisfying
 \eq \label{eq:Palm2}
\E\left(\iint_{\G^2}
f(x,y;\Xi-\d_x-\d_y)\Xi(dx)(\Xi-\d_x)(dy)\right) = \iint_{\G^2} \E
f(x,y;\Xi_{xy}) \lt(dx,dy)
 \en
for any measurable function $f:\G^2\times\H\rightarrow\real_+$, where
$\lt(dx,dy)=\E\Xi(dx)(\Xi-\d_x)(dy)$ is called the {\it second order factorial moment measure} of $\Xi$ [Kallenberg~(1983), \S12.3]. The second order reduced Palm distribution $\L\Xi_{xy}$ can also be viewed as the Radon-Nikodym derivative
$$\prob(\Xi_{xy}\in B)=\frac{\E[\Xi(dx)(\Xi-\d_x)(dy)1_{\{\Xi-\d_x-\d_y\in
B\}}]}{\E\Xi(dx)(\Xi-\d_x)(dy)},\mbox{ for all }B\in\BH.
$$

For $\xi\in\H$ and a Borel set $B\subset \G$, we denote $\xi|_B$ as the restriction of $\xi$
to $B$,
i.e. $\xi|_B(C)=\xi(B\cap C)$ for all Borel sets $C\subset\Gamma$. We call $\{A_x : x\in\G\}$ {\it a type-I neighbourhood} if $x\in A_x\in\BG$ for all $x\in\G$ and
the mapping
 $$
\G\times\H\rightarrow\G\times\H: (x,\xi)\mapsto(x,\xi \vert_{A_x^c})
 $$
is product measurable [see Chen \& Xia~(2004), pp. 2547--2548 for further discussions]. We say that $\{A_{xy} : x,y\in\G\}$ is {\it a type-II neighbourhood} if $\{x,y\}\subset A_{xy}\in\BG$ for all $x,y\in\G$ and
the mapping
 $$
\G^2\times\H\rightarrow\G^2\times\H: ((x,y),\xi)\mapsto((x,y), \xi
\vert_{A_{xy}^c})
 $$
is product measurable. We now define the locally dependent
structures studied in this paper.

\begin{Def} {\rm A point process $\Xi$ is said to satisfy the
{\it type-I local dependence} if there exist two type-I
neighbourhoods $\{A_x : x \in \G \}$ and $\{B_x : x \in\G\}$ such
that $A_x\subset B_x$,
$\L\left(\XixAxc\right)=\L\left(\XiAxc\right)$, $\XiBxc$ is
independent of $\XiAx$, and $\XixBxc$ is independent of
$\XixAx$ for all $x\in\G$. A point process $\Xi$ is said to satisfy
the {\it type-II local dependence} if there exist two type-II
neighbourhoods $\{A_{xy} :  x,y\in\G\}$ and $\{B_{xy}: x,y\in\G\}$
such that $A_{xy}\subset B_{xy}$,
$\L\left(\Xi_{xy}|_{A_{xy}^c}\right)=\L\left(\Xi|_{A_{xy}^c}\right)$,
$\XiBxyc$ is independent of $\Xi|_{A_{xy}}$, and $\XixyBxyc$
is independent of $\XixyAxy$ for all $x,y\in\G$.}
\end{Def}

The locally dependent structures introduced here are parallel
to, but a little stronger than, those in Chen \& Shao~(2004).
The condition $\L\left(\XixAxc\right)=\L\left(\XiAxc\right)$ can be
loosely interpreted as $\Xi(dx)$ is independent of $\XiAxc$. One may
easily establish sufficient conditions for the locally dependent
structures by imposing conditions on neighbourhoods containing balls
[see the descriptive definitions in Barbour \& Xia~(2006)].

To state the error estimates of the PBDP approximation to locally
dependent point processes, we need to introduce the following
notations. Let $\g =\{G_1,\dots,G_k\}\subset\BG$ be a partition of
$\G$, and we choose $t_i\in\G$ such that $\sup_{s \in
G_i}d_0(s,t_i)$ is as small as possible, $i=1,\dots,k$. Note that
$t_i$, regarded as the `designated center' of the set $G_i$, is not necessarily in $G_i$. We define $\m
\eta:=\sum_{i=1}^k\eta(G_i)\d_{t_i}$ for $\eta\in\H$. The mapping is
to `assemble' all the points of
the configuration $\eta$ in each $G_i$ to its
center $t_i$. If we set $d_0(\g)$ as
$$d_0(\g)=\max_{1\le i\le k}\sup_{s\in G_i}d_0(s,t_i),$$
then it is easy to check that
 \eq
d_1(\eta,\m\eta)\le d_0(\g). \label{shuffle1}
 \en
Let $u$ be a positive constant to be chosen in applications, and we take $u=2$ for our examples in Section~\ref{applications}. Let $\F_{TV}$ be  the set of indicator functions of all sets in $\BH$.
For a point process $\Xi$, we define
 \eqs
r_x(\Xi)&:=&4\prob\left.\left(\Xi(B_x^c)+1\le \frac
au\right|\XiBx\right)
 \\ &&
+\frac{4u+10}{a}\max_{1\le j\le
k}\sup_{f\in\F_{TV}}\left|\E\left[\left.f\big(\m(\XiBxc)\big)-f\big(\m(\XiBxc)+\d_{t_j}\big)\right|\XiBx\right]\right|.
 \ens
Similarly, $\bar r_x(\Xi)$ is defined by replacing all the conditional expectations/probability in the definition of $r_x(\Xi)$ with expectations/probability. It is worthwhile to point out that the type-I local dependence implies $\bar r_x(\Xi)=\bar r_x(\Xi_x)$. Let
\eqs
\epsilon_{1,x}(\Xi)&=&r_x(\Xi)\Xi(A_x)\Xi(B_x\setminus A_x)+\bar
r_x(\Xi) \big[ \Xi(A_x)+1 \big] \Xi(A_x)/2+\Xi(A_x)\E
[r_x\left(\Xi\right)\Xi(B_x)],
 \\
 \epsilon_{1,x}(\Xi_x)&=&{r_x(\Xi_x)\Xi_x(A_x)\Xi_x(B_x\setminus A_x)+\bar
r_x(\Xi_x) \big[ \Xi_x(A_x)+1 \big] \Xi_x(A_x)/2+\Xi_x(A_x)\E
[r_x\left(\Xi\right)\Xi(B_x)],}
 \\
\epsilon_{2,x}(\Xi)&=&r_x(\Xi)\Xi(B_x\setminus A_x)+\bar r_x
(\Xi)+\E [r_x(\Xi)\Xi(B_x)],\\
{\epsilon_{2,x}(\Xi_x)}&{=}&{r_x(\Xi_x)\Xi_x(B_x\setminus A_x)+\bar r_x
(\Xi_x)+\E [r_x(\Xi)\Xi(B_x)].}
 \ens
In terms of the type-II local dependence, we define $r_{x,y}$ and $\bar r_{x,y}$ in the same way as $r_x$ and $\bar r_x$ respectively, but with $B_x$ replaced by $B_{xy}$. We then set
 \eqs
\epsilon_{1,x,y}(\Xi)&=&r_{x,y}(\Xi)\Xi(A_{xy})\Xi(B_{xy}\setminus A_{xy})+\bar r_{x,y}(\Xi)(\Xi(A_{xy})+1)\Xi(A_{xy})/2\\
&&+\Xi(A_{xy})\E[ r_{x,y}\left(\Xi\right)\Xi(B_{xy})],\\
{\epsilon_{1,x,y}(\Xi_{x,y})}&{=}&{r_{x,y}(\Xi_{x,y})\Xi_{x,y}(A_{xy})\Xi_{x,y}(B_{xy}\setminus A_{xy})+\bar r_{x,y}(\Xi_{x,y})(\Xi_{x,y}(A_{xy})+1)\Xi_{x,y}(A_{xy})/2}\\
&&{+\Xi_{x,y}(A_{xy})\E[ r_{x,y}\left(\Xi\right)\Xi(B_{xy})],}\\
\epsilon_{2,x,y}(\Xi)&=&r_{x,y}(\Xi)\Xi(B_{xy})+\bar r_{x,y}
(\Xi) + \E \big[ r_{x,y}\left(\Xi\right)\Xi(B_{xy}) \big],\\
{\epsilon_{2,x,y}(\Xi_{x,y})}&{=}&{r_{x,y}(\Xi_{x,y})\Xi_{x,y}(B_{xy})+\bar r_{x,y}
(\Xi_{x,y}) + \E \big[ r_{x,y}\left(\Xi\right)\Xi(B_{xy}) \big].}
 \ens


 \begin{thm} \label{ld1} Assume that the point process $\Xi$ on $\G$ with finite mean measure $\l$ satisfies $\var(|\Xi|)\ge\E|\Xi|$ and the type-I local dependence. Let
$\nu(dx)=\l(dx)/|\l|$, $b=[\var(|\Xi|)-\E|\Xi|]/\var(|\Xi|)$, $a=(1-b)|\l|$, then
$$ d_2(\L\Xi,\bp_{a,b;0;\nu})\le 2d_0(\g)+\int_\G\E \big[ (1+b)(\epsilon_{1,y}(\Xi_y)+\epsilon_{1,y}(\Xi))+b\bar r_y(\Xi)\Xi_y(A_y)+b\epsilon_{2,y}(\Xi_y) \big]\l(dy). $$
 \end{thm}

\begin{thm} \label{ld2} Assume the point process $\Xi$ on $\G$ with finite mean measure $\l$ satisfies $\var(|\Xi|)<\E|\Xi|$, the type-I and type-II local   dependence. Let
\eq\b=\frac{|\l|-\var(|\Xi|)}{|\l|-\var(|\Xi|)+\E|\Xi|^3-(|\l|+1)\E|\Xi|^2},\ a=|\l|+\b(\E|\Xi|^2-|\l|),\label{ld2-1}\en
and
\eq\nu(dx)=\frac1a\left(\l(dx)+\b\int_{y\in\G}\lt(dx,dy)\right).\label{ld2-2}\en
If $\b\ge 0$, then
 \eqs d_2(\L\Xi,\bp_{a,0;\b;\nu})&\le& 2d_0(\g)+\int_\G\E\left(\epsilon_{1,x}(\Xi_x)+\epsilon_{1,x}(\Xi)\right)\l(dx) \\
 &&+\b\iint_{\G^2}\E\left(\epsilon_{1,x,y}(\Xi_{xy})+\epsilon_{1,x,y}(\Xi)+\epsilon_{2,x,y}(\Xi_{xy})\right)\lt(dx,dy).
 \ens
 \end{thm}

\begin{rem} {\rm When one applies these theorems, it is advisable to leave the choice of $\g$ to the last stage so that an optimal bound with the best possible order can be achieved.
}
\end{rem}

A less noticeable fact is that if  one takes $d_0(x,y)=0$, i.e. a
pseudometric on $\G$, and $\g=\{\G\}$, then $d_2$
reduces to $d_{tv}$ for the total counts of point processes,
so our theorems also cover the \PBDP\ approximation to the total
counts of locally dependent point processes in the total variation
distance.

The proofs of the two theorems will be given in sections~\ref{prooftheoremld1} and \ref{prooftheoremld2}. In the next section, let us look at three examples to see how the theorems perform in applications.

\section{Applications}
\label{applications}
\setcounter{equation}{0}

\subsection{Bernoulli process} \label{sectionBernoulli}

Let $\G=[0,1]$, $d_0(x,y)=|x-y|$, $I_1$, \dots, $I_n$ be independent Bernoulli random variables with
$$\prob(I_i=1) = 1-\prob(I_i=0)=p_i,\ 1\le i\le n.$$
Define $\Xi=\sum_{i=1}^nI_i\d_{i/n}$. This simple point process is
particularly useful for proving the Poisson process limit theorems
for the extreme value theory [Embrechts, Kl\"uppelberg and
Mikosch~(1997), Chapter~5]. It was proved in Xia~(1997),
Proposition~3.6 [see also Ruzankin~(2004)], that the accuracy of
Poisson process approximation to $\L\Xi$ is of order
$\sum_{i=1}^np_i^2/\sum_{i=1}^np_i$ and the order can not be
improved when $n$ becomes large. When $p_i$'s are equal to $p$, Xia
\& Zhang~(2008), making use of the symmetric nature of the
distribution $\L\Xi$, proved that an appropriate PBDP can
approximate $\L\Xi$ with approximation error of order
$\left(\frac1n+p\right)\wedge\frac1{\sqrt{np}}$. However, when
$p_i$'s are not the same, the techniques employed in Xia \&
Zhang~(2008) will not work and we demonstrate below that  our
theorems can be applied to this case.

First of all, it is easy to verify that $\Xi$ has mean measure $\l(dx)=\sum_{i=1}^np_i\d_{i/n}(dx)$ and its second order factorial moment measure is $\lt(dx,dy)=\sum_{1\le i\ne j\le n}p_ip_j\d_{i/n}(dx)\d_{j/n}(dy)$. Clearly, $\E|\Xi|>\var(|\Xi|)$, so we can apply Theorem~\ref{ld2}  to estimate the approximation error for $\L\Xi$.

To identify the approximating PBDP distribution, we let
\eqs
\l_j&=&\sum_{i=1}^np_i^j,\ j\ge 2,\\
\b&=&\frac{\l_2}{|\l|^2-\l_2-2|\l|\l_2+2\l_3},\\
a&=&|\l|+\b(|\l|^2-\l_2),
\ens
[cf Brown and Xia~(2001), Theorem~3.1] and
$$\nu(dx)=\frac1a\left(\l(dx)+\b\int_{y\in\G}\lt(dx,dy)\right)=\frac1a\left(\l(dx)+\b\sum_{i=1}^n(|\l|-p_i)p_i\d_{i/n}(dx)\right).$$

Next, we set up an appropriate partition $\g$ of
$\G=\{G_1,\dots,G_k\}$. Let $1\le u_1,\dots,u_k\le n$ such that
$u_1+\dots+u_k=n$, $s_0=0$, $s_j=s_{j-1}+ u_j$ for $1\le j\le k$.
Set $G_1=[0,s_1/n]$ and
$G_j=\left(\frac{s_{j-1}}{n},\frac{s_j}{n}\right]$ for $2\le j\le
k$. We choose $t_j$ as the middle point of the interval $G_j$, $1\le
j\le k$, so that $d_0(\g)=\max_{1\le j\le k} u_j/(2n).$ Define
$W_j=\sum_{i=s_{j-1}+1}^{s_j}I_i$, $1\le j\le k$ and
\eqs\kappa&:=&\max_{1\le j\le k}\max_{s_{j-1}+1\le l_1\ne l_2\le s_j}d_{tv}(\L(W_j-I_{l_1}-I_{l_2}),\L(W_j-I_{l_1}-I_{l_2}+1))\\
&\le&\max_{1\le j\le k}\max_{s_{j-1}+1\le l_1\ne l_2\le s_j}1\wedge\frac{1}{2\sqrt{\sum_{l=s_{j-1}+1}^{s_j}p_l(1-p_l)-p_{l_1}(1-p_{l_1})-p_{l_2}(1-p_{l_2})}},
\ens
where the inequality is due to Lemma~1 of Barbour and Jensen~(1989).
We take $A_x=B_x=\{x\}$, $A_{xy}=B_{xy}=\{x,y\}$, $u=2$, then $\Xi_x(A_x)=\Xi_x(B_x)=\Xi_{xy}(A_{xy})=\Xi_{xy}(B_{xy})=0$,
$$\prob\left(|\Xi|-I_{l_1}-I_{l_2})\le \frac a2\right)\le O(|\l|^{-2}),$$
hence all of $r_x$, $\bar r_x$, $r_{x,y}$ and $\bar r_{x,y}$ are bounded by $O(\kappa/a)$. Applying Theorem~\ref{ld2} gives the following estimate.

\begin{thm} \label{Bern1} With the above setup, if $\b\ge 0$, then
 $$ d_2(\L\Xi,\bp_{a,0;\b;\nu})\le \max_{1\le j\le k} u_j/n+O(\kappa
 \l_2/|\l|). $$
\end{thm}

As a special case, we now assume $p_i$'s are equal to $p$, and take $k=O\left((n(1-p)/p)^{1/3}\right)$, $u_j=O\left((pn^2/(1-p))^{1/3}\right)$, $j=1,\dots,k$, then
$$\kappa=O\left(1\wedge\frac1{(np^2(1-p))^{1/3}}\right).$$
Hence, the following corollary is immediate.

\begin{cor} For the Bernoulli point process $\Xi=\sum_{i=1}^nI_i\d_{i/n}$, where $\{I_i,\ 1\le i\le n\}$ are independent and identically distributed Bernoulli random variables with $\prob(I_1=1)=p$, let $\b=\frac{1}{(n-1)(1-2p)}$, $a=np(1-p)/(1-2p)$, $\nu(dx)=\frac{1}{n}\sum_{i=1}^n\d_{i/n}(dx)$, then
\eq d_2(\L\Xi,\bp_{a,0;\b;\nu})\le O\left(\frac{p^{1/3}}{(n(1-p))^{1/3}}\right),\label{Bern1-2}\en
provided $p<1/2$.
\end{cor}

 \begin{rem}
{\rm The bound \Ref{Bern1-2} is not as good as the bound
$O\left(\left(\frac1n+p\right)\wedge\frac1{\sqrt{np}}\right)$
derived in Xia \& Zhang~(2008) when $p$ is fixed and $n$ becomes
large. Nevertheless, our method does not rely on the specific
symmetric structure of the Bernoulli process $\Xi$ and the bound is
still valid even if the success probabilities for the Bernoulli
random variables vary moderately.}
 \end{rem}

\begin{rem}{\rm A Poisson process approximation to the Bernoulli process is justified when $p\to0$ and $np\to\lambda$. However, in applications of extreme value theory, the value $p$ is often fixed while $n$ is large, so our theory provides a more practical alternative.}
\end{rem}

\subsection{Compound Poisson process}

Barbour \& M{\aa}nsson~(2002) considered compound Poisson process
approximation in $d_2$ distance. The Stein factors for both compound
Poisson random variable and process approximations are generally too
crude to use unless they are sufficiently close to their Poisson
counterparts or satisfy some other restrictive conditions. In this
example, we will show that our PBDP, suitably chosen, will converge
to the compound Poisson process when its cluster distribution is
fixed and the mean of the Poisson process component becomes large,
regardless of whether the compound Poisson process is sufficiently
close to a Poisson process or not.

To begin with, let $\Xi=\sum_{i=1}^\infty iX_i$, where $\{X_i\}$ are
independent Poisson processes on $\G$ with mean measures $\{\mu_i\}$
respectively. For brevity, we write $\Xi\sim{\rm
CP}(\mu_1,\mu_2,\dots)$. It is easy to see that $\var(|\Xi|)\ge
\E|\Xi|$ with equality holds if and only if $\mu_j=0$ for all $j\ge
2$.

Suppose that we have a partition $\g=\{G_1,\dots,G_k\}$ of $\G$.

\begin{thm}\label{compp1}
Let $\l(dx)=\sum_{i=1}^\infty i\mu_i(dx)$, $\nu(dx)=\l(dx)/|\l|$, and
$$b=\frac{\sum_{i=2}^\infty i(i-1)|\mu_i|}{\sum_{i=1}^\infty i^2|\mu_i|},\ a=\frac{|\l|^2}{\sum_{i=1}^\infty i^2|\mu_i|}.$$
Then
\eq d_2({\rm CP}(\mu_1,\mu_2,\dots),\bp_{a,b;0;\nu})\le
O\left(\max_{1\le i\le k}1\wedge \frac{1}{\sqrt{\mu_1(G_i)}}\right)\frac{\sum_{i=1}^\infty i^3|\mu_i|}{a}+2d_0(\g).
\label{compp1-1}\en
\end{thm}

\begin{rem}\label{compp2}
{\rm Suppose the cluster distribution is fixed everywhere and $\mu_1(G)\to\infty$ for every $G\in\BG$ such that $\mu_1(G)>0$, then the upper bound given in \Ref{compp1-1} has the order $o(1)$. To this end, one can partition $\G$ into sets with diameters small enough, then for each set $G_i$ with $\mu_1(G_i)>0$, one can find $\mu_1 (G_i)$ as large as one wishes. Furthermore, suppose $\G$ is a simply connected domain in $\mathbb{R}^d$ with smooth boundary, $d_0(x,y)=|x-y|\wedge 1$, and $\mu_1$ is proportional to the Lebesgue measure, i.e. points are homogeneous on $\G$. Then, the upper bound given in \Ref{compp1-1} has the order $O\left(|\mu_1|^{-\frac1{d+2}}\right)$. As a matter of fact, one can partition $\G$ into boxes with the same diameter of order $O \left(|\mu_1|^{- \frac 1 {d+2}} \right)$, then combine the parts at the boundary of $\G$ to their adjacent boxes totally belonging to $\G$, to obtain ${\cal G}$.}
\end{rem}

\begin{rem}\label{compp3} {\rm The other possible way to define negative binomial process is through a compound Poisson process having a Poisson process of clusters and each cluster carries a random number of points that follows a logarithmic distribution. Remark~\ref{compp2} ensures that if the logarithmic distribution for the clusters is fixed and the Poisson process is homogeneous, then the process will converge to our PBDP distribution when the mean measure of the Poisson process becomes large.}
\end{rem}

\noindent{\it Proof of Theorem~\ref{compp1}.} A measure $\mu$ is called {\it diffuse} if for every point $x\in\G$, $\mu(\{x\})=0$. If $\{\mu_i\}$ are not diffuse, we can enlarge the space $\G$ if necessary and take diffuse measures $\mu_i^n$ such that $|\mu_i^n|=|\mu_i|$ for $i\ge 1$ and $\max_{i\ge 1}d_1(\mu_i^n,\mu_i)\to0$ as $n\to\infty$. We then apply the Kantorovich-Rubinstein duality theorem [Rachev~(1991), Theorem~8.1.1,
p.~168] to couple two sequences of $\G$-valued random elements $\tau_{ij}\sim\mu_i/|\mu_i|$ and $\tau_{ij}^n\sim\mu_i^n/|\mu_i^n|$, $i,j\ge 1$, such that
$$\E d_0(\tau_{ij},\tau_{ij}^n)=d_1\left(\mu_i/|\mu_i|,\mu_i^n/|\mu_i^n|\right)=d_1\left(\mu_i,\mu_i^n\right),$$
and $\{(\tau_{ij},\tau_{ij}^n),\ i,j\ge 1\}$ are independent and independent of $\{X_i,\ i\ge 1\}$. Let $\Xi^n\sim{\rm CP}(\mu_1^n,\mu_2^n,\dots)$, then
\eqs
d_2(\L\Xi,\L(\Xi^n))&\le& \E d_1\left(\sum_{i=1}^\infty i\sum_{j=1}^{|X_i|}\d_{\tau_{ij}},\sum_{i=1}^\infty i\sum_{j=1}^{|X_i|}\d_{\tau_{ij}^n}\right)\\
&\le& \E\left(\frac{\sum_{i=1}^\infty i\sum_{j=1}^{|X_i|}d_0\left(\tau_{ij},\tau_{ij}^n\right)}{\sum_{i=1}^\infty i|X_i|}\right)\\
&\le&\E\left(\frac{\sum_{i=1}^\infty i\sum_{j=1}^{|X_i|}d_1\left(\mu_i,\mu_i^n\right)}{\sum_{i=1}^\infty i|X_i|}\right)\\
&\le& \max_{i\ge 1}d_1(\mu_i^n,\mu_i).
\ens
This observation, together with Lemma~\ref{lem2pbd1}, ensures that we can assume, without loss of generality, that $\{\mu_i\}$ are all diffuse. Otherwise, we can approximate each $\Xi^n$ with a suitable PBDP distribution and then take the limits.

Direct computation gives
\eqs
|\l|&=&\sum_{i=1}^\infty i|\mu_i|,\ \var(|\Xi|)=\sum_{i=1}^\infty i^2|\mu_i|,\\
b&=&\frac{\sum_{i=2}^\infty i(i-1)|\mu_i|}{\sum_{i=1}^\infty i^2|\mu_i|},\ a=\frac{|\l|^2}{\sum_{i=1}^\infty i^2|\mu_i|}.
\ens

Because the compound Poisson process has independent increments, we
let $A_x=B_x=\{x\}$, then
 \eqs
& & r_x(\Xi) =  \bar r_x(\Xi)= r_x(\Xi_x)=\bar r_x(\Xi_x)
 \\ & = &
4\prob\left(|\Xi|+1\le \frac{a}u\right)+\frac{4u+10}{a}\max_{1\le
j\le k}d_{TV}\left(\L(\m\Xi),\L\left( \m\Xi+\d_{t_j}\right)\right),
 \ens
where for any two point process distributions ${\bf P}$ and
${\bf Q}$ on $\H$, $d_{TV} ({\bf P}, {\bf Q}) : = \inf_{\xi \sim
{\bf P}, \eta \sim {\bf Q}} \prob (\xi \neq \eta)$. Noting that $\{
\mu_i \}$ are all diffuse and consequently $\Xi (\{ x \}) = 0$ a.s.
for each $x \in \G$, we have
 \eq
\epsilon_{1,x}(\Xi) = 0, \
\epsilon_{1,x}(\Xi_x)=\bar
r_x(\Xi_x)(\Xi_x(\{x\})+1)\Xi_x(\{x\})/2,\
\epsilon_{2,x}(\Xi_x)=\bar r_x(\Xi_x) . \label{E.CPepsilon} \en

It is well-known that if $Y$ follows Poisson distribution with mean $c$, then
$$d_{tv}(\L Y,\L(Y+1))=\max_{j\ge 0}\prob(Y=j)\le \frac{1}{\sqrt{2ec}},$$
[see Barbour, Holst \& Janson~(1992), Proposition~A.2.7, p.~262].
Hence, we have
 $$
d_{TV}\big(\L(\m\Xi),\L(\m\Xi+\d_{t_j})\big) \le d_{tv}\big(
\L(\Xi(G_j)),\L(\Xi(G_j)+1) \big) \le \frac{1}{\sqrt{2e\mu_1(G_j)}}
\le \frac 1 {\sqrt{\mu_1 (G_j)}} . $$ It is easy to see that we can
write $|\Xi|=\sum_{i=1}^V\eta_i$, where all the random variables $V$
and $\eta_i$'s are independent, $V\sim{\rm Poisson}(|\mu'|)$ with
$\mu'=\sum_{i=1}^\infty\mu_i$, and $\eta_i$'s have the same
distribution $\prob(\eta_i=j)=|\mu_j|/|\mu'|,\ j\ge 1$. If we take
$u=2$, noting that $a\le|\mu'|$, we have
$$\prob\left(|\Xi|+1\le\frac a2\right) \le \prob\left(V \le \frac {|\mu'|} 2 \right) \le O(|\mu'|^{-2}) \le O(a^{-2}).$$
Hence,
\eq
\bar r_x(\Xi_x)=O\left(a^{-1}\max_{1\le i\le k}1\wedge \frac{1}{\sqrt{\mu_1(G_i)}}\right).\label{compp1-2}
\en

Using the independent increments again, we get
\eqs
\var(|\Xi|)&=&\E\int_\G(|\Xi|-|\l|)\Xi(dx)=\int_\G \E (|\Xi_x|+1-|\l|)\l(dx)=\int_\G \E (\Xi_x(\{x\})+1)\l(dx),\\
\E(|\Xi|-1)|\Xi|^2&=&\E\int_\G|\Xi|(|\Xi-\d_x|)\Xi(dx)=\int_\G \E (|\Xi_x|+1)|\Xi_x|\l(dx)\\
&=&|\l|\E|\Xi|^2+2|\l| \int_\G \E \Xi_x(\{x\})\l(dx)+|\l|^2+\int_\G\E(\Xi_x(\{x\})+1)\Xi_x(\{x\})\l(dx),
\ens
which in turn imply
 \eq
\int_\G \E \Xi_x(\{x\})\l(dx)=\var(|\Xi|)-|\l|,\label{compp1-3}
 \en
 \eq
\int_\G\E(\Xi_x(\{x\})+1)\Xi_x(\{x\})\l(dx)=\E(|\Xi|-|\l|)^3-\var(|\Xi|).\label{compp1-4}
 \en

Applying Theorem~\ref{ld1}, (\ref{E.CPepsilon}-\ref{compp1-4}),
together with $0 \le b < 1$, gives \eqa
&&d_2({\rm CP}(\mu_1,\mu_2,\dots),\bp_{a,b;0;\nu})\non\\
&\le&2d_0(\g)+\int_\G \E \left[ \bar r_x (\Xi_x) \left( \big( \Xi_x(\{x\})+1 \big) \Xi_x(\{x\})+\Xi_x(\{x\})+1\right)\right]\l(dx)\non\\
&\le&2d_0(\g)+O\left(a^{-1}\max_{1\le i\le k}1\wedge \frac{1}{\sqrt{\mu_1(G_i)}}\right)\int_\G\E\left[\big(\Xi_x(\{x\})+1\big)\Xi_x(\{x\})+\Xi_x(\{x\})+1\right]\l(dx)\non\\
&=&2d_0(\g)+O\left(a^{-1}\max_{1\le i\le k}1\wedge \frac{1}{\sqrt{\mu_1(G_i)}}\right)\E(|\Xi|-|\l|)^3 . \label{compp1-5}
\ena
Finally, one can verify directly that
$$\E(|\Xi|-|\l|)^3=\E\left(\sum_{i=1}^\infty i(|X_i|-|\mu_i|)\right)^3=\E\sum_{i=1}^\infty i^3(|X_i|-|\mu_i|)^3=\sum_{i=1}^\infty i^3|\mu_i|,$$
which, together with \Ref{compp1-5}, implies \Ref{compp1-1}. \qed

\subsection{Runs}

In the final example, we consider the point process of $k$-runs of 1's in a sequence of independent and identically distributed Bernoulli random variables [cf Example~5.2 of Barbour \& M{\aa}nsson~(2002), p.~1527]. It is easy to see from our derivation that, at the cost of more notational complexity, one can lift the assumption of identical distribution.

To begin with, let $I_1, \cdots, I_n$ be independent Bernoulli random variables
with identical distribution
 $$
\prob(I_i=1) = 1-\prob(I_i=0) = p, \ 1 \le i \le n.
 $$
Let $X_i = \prod_{j=i}^{i+k-1}I_j$ with $k\ge 2$,
where we take $I_j=I_{j-n}$ for $j>n$ to avoid the edge effect. We
define the point process of runs as
 $$\Xi = \sum_{i=1}^n
X_i \d_{i/n}  $$ on $\G=[0,1]$, with 0 being identified the same as
1 and {the distance on the circle $d_0(x,y) = |x-y| \wedge (1 - |x-y|)$}. A point of $\Xi$ at location $i/n$ indicates that there is a
run of $k$ 1's starting at index $i$ and it is clear that the run
may overlap with others around it. Wang \& Xia~(2008) demonstrated
that $\var(|\Xi|)\ge \E |\Xi|$ if and only if
$2+(2k-1)p^k-(2k+1)p^{k-1}\ge 0$, and the latter is easily satisfied
if $p<2/3$. Hence we only consider negative binomial process
approximation to the distribution of $\Xi$.

\begin{thm}\label{runsthm}
Let $k\ge 2$ be a fixed integer,
$$a=\frac{(1-p)np^k}{1+p-(2k+1)p^k+(2k-1)p^{k+1}},\ b=\frac{p\big[2-(2k+1)p^{k-1}+(2k-1)p^k\big]}{1+p-(2k+1)p^k+(2k-1)p^{k+1}},$$
and $\nu(dx)=\frac1n\sum_{i=1}^n\d_{i/n}(dx)$.
Assume $p<2/3$, then
$$ d_2(\L\Xi,\bp_{a,b;0;\nu})\le\left\{ \begin{array}{ll}
O\left(\frac{p^{2/3}}{(np^k)^{1/3}} \right),&{\rm if} \ np^k\ge 1,\\
O(p),& {\rm if} \ np^k<1.
\end{array}\right.
$$
\end{thm}

 \begin{rem}
{\rm The point process of runs in Barbour \& M{\aa}nsson~(2002),
example~5.2, is defined on the carrier space $\G'=[0,n]$ with
{$0$ being identified the same as $n$ and} metric $\tilde
d_0(x,y)=\left(|x-y|p^k\right)\wedge 1$, where $|\cdot|$ is the
distance on the circle. Although $\tilde d_0$ seems to be a natural
choice in the context of compound Poisson process approximation, it
depends on the mean of the process being approximated. An unexpected
effect is, when the parameters vary, it is impossible to judge from
the error estimates whether the approximations become better or
worse. Another defect of the approach in Barbour \&
M{\aa}nsson~(2002) is that a factor $\ln n$ appears inevitably in
the approximation bound, which makes it useless when $n$ becomes
large. In practical applications, $p$ is often fixed while $n$ tends
to be large so that approximate distributions are needed. Our
approximating distribution uses fewer parameters but achieves
approximation bound that decreases when $p$ becomes small and/or $n$
becomes large. }
 \end{rem}

\noindent{\it Proof of Theorem~\ref{runsthm}.} It's easy to verify that the mean measure of $\Xi$ is $\l(dx)=p^k\sum_{i=1}^n\d_{i/n}(dx)$, $\E|\Xi|=|\l|=np^k$ and $\var(|\Xi|)=\frac{np^k}{1-p}\left(1+p-(2k+1)p^k+(2k-1)p^{k+1}\right)$, hence we set
 \eqs
 \nu & = & \frac 1 n\sum_{i=1}^n \d_{i/n},\\
b & = &\frac{\var(|\Xi|)-\E|\Xi|}{\var(|\Xi|)}=\frac{p\big[2-(2k+1)p^{k-1}+(2k-1)p^k\big]}{1+p-(2k+1)p^k+(2k-1)p^{k+1}},
 \\
a & = & (1-b)np^k=\frac{(1-p)np^k}{1+p-(2k+1)p^k+(2k-1)p^{k+1}}.
 \ens

We assume $|\l|\ge1$ first. To tackle the dependence resulting from
the overlapping runs, we introduce the neighbourhoods
$A_{i/n}=\{j/n:\ i-k+1\le j\le i+k-1\}$ and $B_{i/n}=\{j/n:\
i-2k+2\le j\le i+2k-2\}$, where $j$ is interpreted as $j+n$ if
$j<0$, and $j-n$ if $j>n$. Next, we choose $\g=\{G_j:\ 1\le j\le
l_n\}$ by taking $l_n=O\left(n^{1/3}p^{(k-2)/3}\right)$,
$G_j=(s_{j-1}/n,s_j/n]$ for $j=1,\ \dots, l_n$, where $s_0=0$,
$s_j=s_{j-1}+u_j$ for $j=1,\dots,l_n$ with $u_1,\ \dots,\
u_{l_n}=O\left(n^{2/3}p^{(2-k)/3}\right)$ and
$\sum_{j=1}^{l_n}u_j=n$.

To estimate $r_x$, we take $u=2$, write $x=i/n$ and
$Y_i=\sum_{|j-i|>4k-4}X_j$. Applying
the Bienaym\'e-Chebyshev inequality gives
\eqa
&&\prob\left.\left(\Xi\left(B_x^c\right)+1\le \frac au\right\vert\XiBx\right)\le\prob\left(Y_i+1\le\frac au\right)\le \prob\left(Y_i-\E Y_i\le \frac au-\E Y_i\right)\nonumber\\
&&\le \prob\left(\left|Y_i-\E Y_i\right|\ge \left|\frac a2-(n-(8k-7))p^k\right|\right)
\le \frac{\E(Y_i-\E Y_i)^4}{((1+b)np^k/2-(8k-7)p^k)^4}.\label{runs1-2}
\ena
However,
$$\E(Y_i-\E Y_i)^4=\sum_{|j_v-i|>4k-4,\ v=1,2,3,4}\E\prod_{v=1}^4\left(X_{j_v}-\E X_{j_v}\right)$$
and the summand reduces to 0 if one of the $j_v$'s is not in the neighbourhoods of the others, hence
$$\E(Y_i-\E Y_i)^4\le 12|\l|^2\left(\frac{1-p^k}{1-p}\right)^2+9|\l|(12k-9)=O\left(|\l|^2\right).$$
This, together with \Ref{runs1-2}, implies \eq
\prob\left.\left(\Xi\left(B_x^c\right)+1\le \frac
au\right\vert\XiBx\right)\le
O\left(|\l|^{-2}\right).\label{runs1-3}\en The same argument also
leads to \eq \prob\left.\left(\Xi_x\left(B_x^c\right)+1\le \frac
au\right\vert\XixBx\right)\le
O\left(|\l|^{-2}\right).\label{runs1-4}\en

For $f\in\F_{TV}$, we will show that
\eqa
&& \left|\E\left[f\left(\m\left(\left.\XiBxc\right)\right)-f\left(\m\left(\XiBxc\right)+\d_{t_j}\right)\right|\XiBx\right]\right|
\le O\left(n^{-1/3}p^{-(k+1)/3}\right),\label{runs1-5}\\
&&\left|\E\left[f\left(\m\left(\XiBxc\right)\right)-f\left(\m\left(\XiBxc\right)+\d_{t_j}\right)\right]\right|
\le O\left(n^{-1/3}p^{-(k+1)/3}\right).\label{runs1-6}
\ena
In fact, if we write $t_j=(s_{j-1}+s_j)/(2n)$, $x=i/n$, then there are two cases to consider.

Case 1. $s_{j-1}<i\le s_j$. Because of the symmetry of our argument, we assume without loss of generality that
$i\le \frac{s_{j-1}+s_j}2$. We write ${\bf I}_1=(I_1,\dots,I_{i+2k-2},I_{s_j-k+1},\dots,I_n)$, ${\bf I}_2 = (I_{i+2k-1}, \cdots, I_{s_j-k} )$, ${\bf v}=(v_1,\dots,v_{i+2k-2},v_{s_j-k+1},\dots,v_n)$. For any vector ${\bf v}$ with $v_l\in\{0,1\}$, $\forall l$,
due to Wang and Xia~(2008, Lemma~2.1), the number $W({\bf v},{\bf I}_2)$ of $k$-runs of the sequence
$$v_{s_{j-1}+1},\dots,v_{i+2k-2},I_{i+2k-1},\dots,I_{s_j-k},v_{s_j-k+1},\dots,v_{s_j}$$ satisfies \eq d_{tv}(\L W({\bf v},{\bf I}_2),\L (W({\bf v},{\bf I}_2)+1))\le O\left(n^{-1/3}p^{-(k+1)/3}\right) . \label{runs1-7}\en
For ease of notation, we use $\Xi({\bf v},{\bf I}_2)$  to stand for the point process of $k$-runs of the sequence
$$v_1,\dots,v_{i+2k-2},I_{i+2k-1},\dots,I_{s_j-k},v_{s_j-k+1},\dots,v_n.$$
Then, for $f\in\F_{TV}$,
\eqs
&&\left|\E\left[\left.f\left(\m\left(\XiBxc\right)\right)-f\left(\m\left(\XiBxc\right)+\d_{t_j}\right)\right|{\bf I}_1 = {\bf v}\right]\right|\nonumber\\
&\le & d_{TV} \left( \L \big( \m\left(\Xi\left({\bf v},{\bf I}_2\right)\vert_{B_x^c}\right) \big) , \L \big( \m\left(\Xi\left({\bf v},{\bf I}_2\right)\vert_{B_x^c}\right)+\d_{t_j} \big) \right)\nonumber\\
&=& d_{tv} \left( \L  W({\bf v},{\bf I}_2)  , \L \big( W({\bf
v},{\bf I}_2)+1 \big) \right) , \ens and this, together with
(\ref{runs1-7}), yields that
 $$
\left|\E\left[\left.f\left(\m\left(\XiBxc\right)\right)-f\left(\m\left(\XiBxc\right)+\d_{t_j}\right)\right|\XiBx\right]\right|
\le  O\left(n^{-1/3}p^{-(k+1)/3}\right) .
 $$

Case 2. $i\not\in\left(s_{j-1},s_j\right]$. The proof is omitted since it is essentially the same as {that of} case 1 with some minor change of notations only.

The proof of \Ref{runs1-6} is similar. Now, combining (\ref{runs1-3}-\ref{runs1-6}) yields {$r_x(\Xi)$ and $\bar r_x(\Xi)$ are both bounded by $O\left(|\l|^{-1}\left(n^{-1/3}p^{-(k+1)/3}\right)\right)$.} These, together with some crude estimates, e.g. $\E\Xi_x(A_x)\le
(2k-2)p, \E\Xi_x(A_x)\Xi_x(B_x\backslash A_x)\le (2k-2)^2p $ etc.,
imply {that $\E\epsilon_{1,x}(\Xi_x)$, $\E\epsilon_{1,x}(\Xi)$ and $b\E\epsilon_{2,x}(\Xi_x)$ are all bounded by $O\left(|\l|^{-1}\left(n^{-1/3}p^{-(k+1)/3}\right)\right)p$.}
Therefore, if $|\l|\ge1$, the proof is completed by substituting these estimates for the corresponding terms in Theorem~\ref{ld1}.

Finally, if $|\l|< 1$, we take $l_n=O\left(p^{-1}\right)$, $u_1,\dots,u_{l_n}=O\left(np\right)$. Then the right hand side of \Ref{runs1-3} and \Ref{runs1-4} can be replaced with 0, and the upper bounds for \Ref{runs1-5} and \Ref{runs1-6} become $O(1)$, which in turn imply {that $r_x(\Xi)$ and $\bar r_x(\Xi)$ are both bounded by $O\left(|\l|^{-1}\right)$.} Consequently,
$\E\epsilon_{1,x}(\Xi_x)$, $\E\epsilon_{1,x}(\Xi)$ and
$b\E\epsilon_{2,x}(\Xi_x)$ are all bounded by $O\left(|\l|^{-1}\right)p$. We then employ Theorem~\ref{ld1} to obtain the bound $p$, as claimed.
\qed


\section{Proof of Theorem~\ref{mainresult}.}
\label{proofofsteinfactor}
\setcounter{equation}{0}


The proof of Theorem~\ref{mainresult} relies on the coupling and
analysis techniques. The main obstacle in coupling various
birth-death systems
together is the difficulty of identifying the individual particles
from their locations. To circumvent the repeats of points, we need
to lift the space to a higher-dimensional carrier space and tackle
the problem in the lifted space. Such technique has been proved very
effective in handling this type of situations [Chen and Xia~(2004) and
Xia~(2005)].

\subsection{Lifting the carrier space}

In this subsection, we define $\tG=\G\times[0,1]$ and the pseudometric $\t d_0$ on $\tG$ as
$$\t d_0((x_1,t_1),(x_2,t_2))=d_0(x_1,x_2).$$
Let $\tH$ be the class of all finite integer-valued measures
on $\tG$ and $\t d_1$ be the induced pseudometric from $\t d_0$ in the same way as $d_1$
from $d_0$. For $\t\xi\in\tH=\sum_{i=1}^n\d_{(x_i,t_i)}$, we define
$\t\xi\vert_\G=\sum_{i=1}^n\d_{x_i}$ and extend a function $f\in\F$ to a function on $\tH$ by
$$\t f(\t \xi)=f(\t\xi\vert_\G).$$
It is not hard to check that for each $f\in\F$, $\t f$ is a $\t d_1$-Lipschitz function:
$|\t f(\t\xi_1)-\t f(\t\xi_2)|\le \t d_1(\t\xi_1,\t\xi_2)$ for all $\t\xi_1,\t\xi_2\in\tH$.

Next, we define $\t\mu$ as the product measure of $\mu$ and Lebesgue measure on
$[0,1]$. Regardless of whether $\mu$ is diffuse,
the measure $\t\mu$ is always diffuse on $\tG$. Let
 \eqs
\tA\t h(\t\xi)&=&(a+b|\t\xi|)\int_{\tG}(\t h(\t\xi+\d_{\t x})-\t h(\t\xi))\t\mu(d\t x)\\
&&+(1+\beta(|\t\xi|-1))\int_{\tG}(\t h(\t\xi-\d_{\t x})-\t h(\t
\xi))\t\xi(d\t x).
 \ens
Birth-death systems on $\tG$ with the generator $\tA$ evolve in
the same way as birth-death systems on $\G$ with the
generator $\A$.

To carry out the proof of Theorem~\ref{mainresult}, for a given
birth-death system $\Z_{\xi} (\cdot)$ with
$\xi=\sum_{i=1}^n\d_{x_i}$, one can lift it to $\t\Z_{\t\xi} (\cdot)$ by setting a $\t\xi\in\tH$ consisting of distinct
particles at $(x_i,t_i)$, $1\le i\le n$, where $t_1$, $\dots$, $t_n$
are distinct elements of $[0,1]$, and throwing each new born
particle at $z$ equally likely onto $\{z\}\times[0,1]$,
independently of the others. Then,
$$\t f(\t\Z_{\t\xi} (t))=f(\t\Z_{\t\xi} (t) \vert_\G )=f(\Z_{\xi}(t)), \ \ \ \forall \ t \ge 0.$$
This procedure enables us to assume
from now on that, without loss of generality, $\mu$ is diffuse and
the particles at $\xi$, $\eta$, $x$ and $y$ are all distinct.

\subsection{Proof of \Ref{estimateCn}}

First of all, the proportion of the surviving initial particles at time $t$ can be estimated as
 \begin{equation}
\E \frac{|\eta \cap \Z_\eta (t)|}{| \Z_\eta (t) |} \le \min \left\{
\( 1 + \frac{a}{2 |\eta|}  (e^t - 1) \) ^{-1} , e^{- (a \wedge b) t}
\right\} .\label{ration}
 \end{equation}

To this end, we define $ g (\zeta) := |\eta \cap \zeta| /
|\zeta| $ for the fixed $\eta \in \H$, where $0/0$ is
interpreted as $0$. Recall that $V(\zeta)$ has the uniform distribution on
the sites in $\zeta$, and we have $ \E g (\zeta - \d_{V(\zeta)}) = g
(\zeta)$. Hence
 $$
\A g (\zeta) = \big( a + b | \zeta | \big) \( \E g ( \zeta + \d_U) -
g (\zeta)  \)
 $$
since the last term in (\ref{generator}) vanishes. {Noticing} that with
probability 1, $U \notin \eta$, we have
 $$
g ( \zeta + \d_U) - g (\zeta) = |\eta \cap \z| \( \frac 1 { |\zeta +
\d_U| } - \frac 1 { |\zeta| }  \) = - \frac{ |\eta \cap \zeta|
}{|\zeta| ( | \zeta | + 1)} \ \ \ {\rm a.s.}
 $$
It follows that
 $$
\A g (\zeta) \le \min \left\{ -  \frac{ a |\eta \cap \zeta|  }{2
|\zeta|^2 } , - (a \wedge b ) g (\zeta) \right\}.
 $$
Therefore, setting $\varphi (t) = \E g(\Z_\eta (t))$, we have
  \begin{equation} \label{neweq1}
\varphi^\prime (t) = \E \A g (\Z_\eta (t)) \le \min \left\{ - \frac
a 2 \E \frac{|\eta \cap \Z_\eta (t) | }{ |\Z_\eta (t) |^2} , - (a
\wedge b) \varphi (t)  \right\} .
 \end{equation}

By the Cauchy inequality,
 $$
\( \E \frac{|\eta \cap \Z_\eta (t)|}{| \Z_\eta (t) |} \) ^2 \le  \E
\frac{|\eta \cap \Z_\eta (t)|}{| \Z_\eta (t) |^2}  \E |\eta \cap
\Z_\eta (t)| \le  |\eta| e^{-t} \E \frac{|\eta \cap \Z_\eta (t)|}{|
\Z_\eta (t) |^2}  ,
 $$
where the second inequality holds since each particle dies with rate
at least 1. Therefore,
 $$
\E \frac{|\eta \cap \Z_\eta (t)|}{| \Z_\eta (t) |^2} \ge  \frac
{e^t} {|\eta|} \( \E \frac{|\eta \cap \Z_\eta (t)|}{| \Z_\eta (t) |}
\) ^2 .
 $$
This, together with (\ref{neweq1}), yields
 $$
\varphi^\prime (t) \le \min \left\{ - \frac {a e^t} {2 |\eta|}
\varphi (t) ^2 , - (a \wedge b) \varphi (t) \right\} .
 $$
Therefore, \Ref{ration} follows from the fact that $\varphi (0) = 1$.

Next, suppose $\eta \in \H$,
$|\eta| = n$ and the particles at $x, y$ and $\eta$ are all distinct. We start with $\Z_{\eta +
\d_x} {(\cdot)}$ and construct $\Z_{\eta + \d_y}  {(\cdot)}$ by replacing $x$ with $y$.
Let $\tau_z  = \inf \{ t : z \notin \Z_{\eta+\d_x} (t)  \} $ for $z\in \eta+\d_x$.
Then, $\Z_{\eta + \d_x} (t) = \Z_{\eta + \d_y} (t)$ for $t \ge
\tau_x$. For $t < \tau_x$,
 $$
|f ( \Z_{\eta + \d_x} (t)) - f( \Z_{\eta+\d_y} (t)) |\le d_1 (\Z_{\eta
+ \d_x} (t), \Z_{\eta + \d_y} (t)) \le \frac 1 {|\Z_{\eta + \d_x}
(t)|} .
 $$
Therefore,
 \begin{equation} \label{deltah}
| h(\eta + \d_x) - h(\eta + \d_y) | \le \int_0^\infty \E \frac{1_{
\{ \tau_x > t \} }} { | \Z_{\eta + \d_x} (t)| } d t .
 \end{equation}
Notice that for all $z \in \eta + \d_x$,
 $$
\E \frac{1_{ \{ \tau_x > t \} }} { | \Z_{\eta + \d_x} (t)| }  = \E
\frac{1_{ \{ \tau_z > t \} }} { | \Z_{\eta + \d_x} (t)| } ,
 $$
 which implies that
 $$
\E \frac{1_{ \{ \tau_x > t \} }} { | \Z_{\eta + \d_x} (t)| }  =
\frac 1 {n+1} \E \frac{ \sum_{z \in \eta + \d_x} 1_{ \{ \tau_z > t
\} }} { | \Z_{\eta + \d_x} (t)| } = \frac 1 {n+1} \E \frac{ | ( \eta
+ \d_x ) \cap  \Z_{\eta + \d_x} (t) | } { | \Z_{\eta + \d_x} (t)| }
.
 $$
This, together with (\ref{deltah}) and {\Ref{ration}}, implies
that
 \begin{eqnarray*}
C_n
 & \le &
\frac 1 {n+1} \int_0^\infty  \frac 1 {1 + \frac{a}{2(n+1)} (e^t - 1)
} d t = \frac{\ln (n+1) - \ln \frac{a}{2}}{n+1 - \frac{a}{2}}
 \\ & \le &
\frac 1 2 \( \frac 1 {n+1} + \frac 2 a \) = \frac{1}{2 (n+1)} +
\frac 1 a ,
 \end{eqnarray*}
where the result also includes the case $a = 2 (n+1)$, and
 $$
C_n \le \int_0^\infty \frac{1}{n+1}   e^{- (a \wedge b) t} d t =
\frac{1}{ (a \wedge b) (n+1)} .
 $$
On the other hand,
 $$
\E \frac{1_{ \{ \tau_x > t \} }} { | \Z_{\eta + \d_x} (t)| } \le \prob
(\tau_x > t) \le e^{-t},
 $$
hence  $C_n \le 1$. \qed

\subsection{Proof of \Ref{mainresult2}}

Suppose $|\xi| = n$ and particles at $\xi$, $\eta$, $x$ and $y$ are all distinct. Recall that $\A h (\xi + \d_x) = f(\xi + \d_x)
- \bp (f)$, i.e.
 \eqa
 & &
\a_{n+1}  \E h(\xi + \d_x+\d_U)  + \b_{n+1} \E h(\xi + \d_x -
\d_{V(\xi + \d_x)}) - \big( \a_{n+1}+ \b_{n+1} \big)  h(\xi+\d_x)
  \non\\ & =&
f (\xi + \d_x) - \bp (f) .\label{xia4.1}
   \ena
It follows that
  \begin{eqnarray*}
 & &
\E h(\xi + \d_x+\d_U)
 \\ & = &
\frac{f(\xi + \d_x) - \bp (f) }{\a_{n+1}} + \frac{\a_{n+1}+
\b_{n+1}} {\a_{n+1}} h(\xi+\d_x)  -  \frac{\b_{n+1}}{\a_{n+1}} \E
h(\xi + \d_x - \d_{V(\xi + \d_x)}) .
   \end{eqnarray*}
Hence
\eqa
\D_2h (\xi;x,y) & = &
h(\xi+\d_x+\d_y)  - \E h(\xi+\d_x+\d_U) + h(\xi +\d_x) - h(\xi+\d_y)
 \non\\ & &
+  \E h(\xi+\d_x+\d_U) - 2 h(\xi+\d_x)  + h(\xi)
 \non\\  & = &
h(\xi+\d_x+\d_y)  - \E h(\xi+\d_x+\d_U) + h(\xi +\d_x) - h(\xi+\d_y)
 \non\\ & &
+ \frac{f(\xi + \d_x) - \bp (f) }{\a_{n+1}} +   \E  \( h(\xi) -  h
\( \xi + \d_x - \d_{V(\xi + \d_x)} \) \)
 \non\\ & &
+ \frac{\a_{n+1}-\b_{n+1}} {\a_{n+1}} \E \(  h \( \xi + \d_x -
\d_{V(\xi + \d_x)} \) - h(\xi+\d_x)  \) .\label{xia4.2}
 \ena
Swapping $x$ and $y$, we get
\eqa
\D_2h (\xi;y,x)
 & = &
h(\xi+\d_y+\d_x)  - \E h(\xi+\d_y+\d_U) + h(\xi +\d_y) - h(\xi+\d_x)
 \non\\ & &
+ \frac{f(\xi + \d_y) - \bp (f) }{\a_{n+1}} +   \E  \( h(\xi) -  h
\( \xi + \d_y - \d_{V(\xi + \d_y)} \) \)
 \non\\ & &
+ \frac{\a_{n+1}-\b_{n+1}} {\a_{n+1}} \E \(  h \( \xi + \d_y -
\d_{V(\xi + \d_y)} \) - h(\xi+\d_y)  \) .\label{xia4.3}
\ena
Since $ \D_2h (\xi;x,y)= \D_2h (\xi;y,x)$ and $|f - \bp (f)| \le 1$, we take the average of \Ref{xia4.2} and \Ref{xia4.3} to reach the bound
 \begin{equation}  \label{1}
| \D_2h (\xi;x,y) | \le \frac{1}{\a_{n+1}} + C_{n-1} + C_{n+1} +
\left| \frac{\a_{n+1}-\b_{n+1}} {\a_{n+1}}  \right| \D_n,
 \end{equation}
where
 $$
\D_n = \sup \{ |h(\eta + \d_x) - h(\eta)| : |\eta| = n, x \in \G \} .
 $$

On the other hand, we use \Ref{xia4.1} again to obtain
 \eqs
 & &
\E h(\xi + \d_x - \d_{V(\xi + \d_x)})
 \\ & = &
\frac{f (\xi + \d_x) - \bp (f)}{\b_{n+1}} + \frac{\a_{n+1}+
\b_{n+1}} {\b_{n+1}}  h(\xi+\d_x) - \frac{\a_{n+1}} {\b_{n+1}} \E
h(\xi + \d_x+\d_U) ,
   \ens
and argue in the same way as for \Ref{1} to get
  \begin{equation} \label{2}
| \D_2h (\xi;x,y)| \le \frac{1}{\b_{n+1}} + C_{n-1} + C_{n+1}  +
\left| \frac{\b_{n+1} - \a_{n+1}} {\b_{n+1}}  \right| \D_{n+1}.
 \end{equation}
In subsection \ref{sectionDh} below, we will prove that
\eq \label{propestimate}
 \left\{ \begin{array}{ll}
\frac{\a_{n+1}-\b_{n+1}} {\a_{n+1}} \D_n \le \frac1{\a_{n+1}}+C_n , & {\rm if \
} \a_{n+1} \ge \b_{n+1}, \\
\frac{\b_{n+1} - \a_{n+1}} {\b_{n+1}} \D_{n+1} \le \frac1{\b_{n+1}}+C_n, & {\rm
otherwise,}
 \end{array} \right.
 \en
and so it follows from \Ref{1} and \Ref{2}
that
 \begin{equation} \label{E.D2h}
| \D_2 h (\xi) | \le C_{n-1} + C_n + C_{n+1} + 2\left( \frac 1
{\a_{n+1}} \wedge \frac 1 {\b_{n+1}} \right) .
 \end{equation}

For $n = 0$, \Ref{estimateCn} enables us to conclude that $C_k \le 1$, $C_{-1} = 0$, and it follows from
\Ref{E.D2h} that
 $$
| \D_2 h (\xi) | \le 2 + \frac 2  a \le \frac 2 {n+1} + \frac 5 a .
 $$
For $n \ge 1$, using the estimate $C_k \le  \frac 1 {2(k+1)} + \frac 1 a$
in \Ref{estimateCn}, the fact $2n \ge n+1$, and the bound given in \Ref{E.D2h}, we have
 $$
| \D_2 h (\xi) | \le \frac 1 {2n} + \frac 1 {2(n+1)} + \frac 1
{2(n+2)} + \frac 5 a \le \frac 2 {n+1} + \frac 5 a.
 $$
This completes the proof of \Ref{mainresult2}. \qed

\subsection{Proof of \Ref{propestimate}} \label{sectionDh}

Since $\{|\Z_\eta (t)|,\ t\ge 0\}$ is a birth-death process with birth rates $\{\a_k\}$, death rates
$\{\b_k\}$ and initial value $|\eta|$, we follow the convention in Brown and Xia~(2001) to define
$\tau_{|\eta|, k} = \inf \{  t: |\Z_\eta (t)| = k \}$, $\tau_m^+
= \tau_{m,m+1}$ and $\tau_m^- = \tau_{m,m-1}$.

For any $\eta \in \H$ with $|\eta| = n$, by the strong Markov property of $\{\Z_\eta (t),\ t\ge 0\}$,
 $$
h(\eta)  =  - \E \int_0^{\tau_n^+} (f (\Z_{\eta} (t) )  - \bp (f) )
d t + \E h(\Z_{\eta} (\tau_n^+) ) ,
 $$
which implies that
 \begin{equation}\label{need1}
\left|  h \( \eta \) - \E h \( \Z_\eta (\tau_n^+) \) \right| \le \E
\tau_n^+ .
 \end{equation}
Now we compare $\E h \( \Z_\eta ( \tau_n^+ ) \) $ with $h (\eta +
\d_x)$. Let $K_n^+$ be the number of particles in $\eta$ that have
died before $\tau_n^+$. Clearly, $0 \le K_n^+ \le n$. Given $K_n^+ = k$,
there are at most $k+1$ pairs of mismatched points between $\Z_{\eta} (\tau_n^+)$ and $\eta + \d_x$, consequently,
 $$
\left|  \E \( h \( \Z_{\eta} (\tau_n^+) \) | K_n^+ = k \) - h(\eta +
\d_x) \right| \le C_n (k+1).
 $$
This in turn leads to
 \eq
\left| \E h \( \Z_{\eta} (\tau_n^+) \) - h(\eta + \d_x) \right| \le
C_n ( \E K_n^+ +1) .\label{xia4.4}
 \en
Combining \Ref{need1} and \Ref{xia4.4} gives
 $$
\D_n \le \E \tau_n^+ + C_n ( \E K_n^+ +1) .
 $$

Likewise, for $\eta \in \H$ with $|\eta| = n+1$, it follows from the strong Markov property of $\{\Z_{\eta+\d_x}(t),\ t\ge 0\}$ that
$$
h(\eta+\d_x)  =  - \E \int_0^{\tau_{n+2}^-} (f (\Z_{\eta} (t) )  - \bp (f) )
d t + \E h(\Z_{\eta+\d_x} (\tau_{n+2}^-) ) ,
 $$
giving
  \eq
| h(\eta + \d_x)  - \E h(\Z_{\eta +\d_x} (\tau_{n+2}^-) ) | \le \E
\tau_{n+2}^- .\label{xia4.5}
 \en
Let $K_{n+2}^- $ be the number of particles in $\eta+\d_x$ that have
died before $\tau_{n+2}^-$, then there are at most $K_{n+2}^- $ mismatched pairs of points between $\Z_{\eta + \d_x} (\tau_{n+2}^-)$ and $\eta$, leading to the bound
 \eq
| \E h( \Z_{\eta + \d_x} (\tau_{n+2}^-)  ) - h(\eta) | \le C_n  \E
K_{n+2}^- .\label{xia4.6}
 \en
Collecting the estimates \Ref{xia4.5} and \Ref{xia4.6}, we obtain
 $$
\D_{n+1} \le  \E \tau_{n+2}^- + C_n \E K_{n+2}^- .
 $$

Put $F(k) = \sum_{i=0}^k \pi_i$ and $\bar{F} (k) = \sum_{i=k}^\infty
\pi_i$. By Lemma~2.2 and Lemma~2.4. in Brown \& Xia~(2001),
 \eqa
&&\E \tau_k^+ = \frac{F(k)} {\a_k \pi_k}, \ \E \tau_k^- =
\frac{\bar{F} (k)}{\beta_k \pi_k}, \non\\
&&\ \frac{F(k)}{F(k-1)} \ge
\frac{\a_k}{\b_k} \ge \frac{\bar{F} (k+1)}{\bar{F} (k)}\label{xia4.7}
 \ena
since $\a_k - \a_{k-1} \le \b_k - \b_{k-1}$ 
for all $k$. It follows from the first inequality of \Ref{xia4.7} that
$$\frac{(\a_{n+1}-\b_{n+1})F(n)}{\a_n\pi_n}\le\frac{\b_{n+1}F(n+1)-\b_{n+1}F(n)}{\a_n\pi_n}=\frac{\b_{n+1}\pi_{n+1}}{\a_n\pi_n}=1,$$
which in turn yields
$$\E \tau_n^+ = \frac{F(n)}{\a_n \pi_n} \le \frac 1 {
\a_{n+1} - \b_{n+1}}, \ \ {\rm if} \ \b_{n+1} < \a_{n+1}.$$
Likewise, using the second inequality of \Ref{xia4.7}, we get
$$\E \tau_{n+2}^-  = \frac{\bar{F} (n+2)}{\beta_{n+2} \pi_{n+2}} \le
\frac 1 {\b_{n+1} - \a_{n+1}}, \ \ {\rm if} \ \b_{n+1} > \a_{n+1}.
$$

To complete the proof of \Ref{propestimate}, it remains to show
 \begin{eqnarray}
\frac{\a_{n+1}-\b_{n+1}} {\a_{n+1}} ( \E K_n^+ +1) \le 1, & {\rm if}
\ \b_{n+1} < \a_{n+1}, \label{need3}
 \\
\frac{\b_{n+1}-\a_{n+1}} {\b_{n+1}} \E K_{n+2}^- \le 1, & {\rm if} \
\b_{n+1} > \a_{n+1}. \label{need33}
 \end{eqnarray}
To this end, we derive a recursive formula for $\E K_m^+$ and $\E
K_m^-$, $m\ge 1$, in Lemma~\ref{EK} later and give their estimates in following
Lemma~\ref{xia4.8}. In particular, since $\a_k - \b_k$ decreases in
$k$ and $\a_k$ increases in $k$, it follows from Lemma~\ref{xia4.8}
that, if $\a_{n+1} > \b_{n+1}$,
 $$
1 + \E K_n^+ \le \frac {\a_n} {\a_n - \b_n} \le \frac {\a_{n+1}}
{\a_{n+1} - \b_{n+1}},
 $$
which is equivalent to (\ref{need3}). On the other hand, noting that
$$\b_{n+2} / (\b_{n+2} -
\a_{n+2}) \le \b_{n+1} / (\b_{n+1} - \a_{n+1}) $$ as $\b_{n+1} -
\a_{n+1} > 0$, applying Lemma~\ref{xia4.8} again, we obtain $\E K_{n+2}^- \le
\b_{n+1} / (\b_{n+1} - \a_{n+1})$ and hence (\ref{need33}) follows. \qed

  \begin{lem} \label{EK}
The following recursive formulae hold for $m \ge 1$:
 $$
\E K_m^+ =  \frac{m \b_m \big( 1 + \E K_{m-1}^+ \big) }{ m \a_m +
\b_m \big( 1 + \E K_{m-1}^+ \big) } , \ \ \ \E K_m^- = 1 +
\frac{(m-1) \a_m \E K_{m+1}^-}{\a_m \E K_{m+1}^- + (m+1) \b_m } .
 $$
 \end{lem}

\noindent{\it Proof.} Noting that all particles die equally likely, an initial particle in the initial configuration $\eta$ with
$|\eta| = m$ dies before $\tau_m^+$
with probability $\frac 1 m \E K_m^+$, and if survives, it dies
before $\tau_{m,m+2}$ with probability $\frac 1 {m+1} \E K_{m+1}^+$.
That is, the probability that an initial particle dies before $\tau_{m,m+2}$ is
 $$
\frac 1 m \E K_m^+ + \( 1 - \frac 1 m \E K_m^+ \)  \frac 1 {m+1} \E
K_{m+1}^+ .
 $$
Therefore, there are in average $ \E
K_m^+ + \frac {m - \E K_m^+} {m+1} \E K_{m+1}^+$ initial particles
die before $\tau_{m,m+2}$.

On the other hand, $K_m^+ = 0$ means that the first change of the configuration of the birth-death system $\Z_\eta (\cdot)$ is a birth, so
$K_m^+ = 0$ with probability $\frac{\a_m}{\a_m + \b_m}$. However, if
the first change is a death, which happens with
probability $\frac{\b_m}{\a_m + \b_m}$, then one particle at some site $x$ of $\eta$ will die at $\btau_\eta=\inf \{t : \Z_\eta (t)
\neq \eta \}$. In the latter case, using the conclusion in the
preceding paragraph, the mean number of particles in $\eta - \d_x$
dying before the birth-death system $\Z_{\eta-\d_x}
(\cdot)$ reaches the size $m+1$ is $\E K_{m-1}^+ + \frac
{m-1 - \E K_{m-1}^+} {m} \E K_{m}^+$. In summary, we have
established the relationship
 $$
\E K_m^+ = \frac{\b_m}{\a_m + \b_m} \( 1 + \E K_{m-1}^+ + \frac {m-1
- \E K_{m-1}^+} {m} \E K_{m}^+ \) ,
 $$
 which is equivalent to the first recursive formula.

The same argument can be adapted to prove the second recursive formula. In fact, assume $|{\eta}| = k \ge 2$, an initial particle
in {$\eta$} dies before $\tau_{k,k-2}$ with probability
 $$
\frac 1 k \E K_k^- + \( 1 - \frac 1 k \E K_k^- \)  \frac 1 {k-1} \E
K_{k-1}^-.
 $$
Now, let $|\eta|=m$. With probability $\frac{\b_m}{\a_m + \b_m}$, the first change of $\Z_\eta {(\cdot)}$ is a death, giving $K_m^-=1$. Assume next that the first change is a birth, then, as shown above, each initial particle dies  before the size reaches $m-1$ with probability $
\frac 1 {m+1} \E K_{m+1}^- + \( 1 - \frac 1 {m+1} \E K_{m+1}^- \)  \frac 1 m \E
K_m^-$. It then follows that
 $$
\E K_m^- = \frac{\b_m}{\a_m + \b_m} + \frac{\a_m}{\a_m + \b_m} \frac
m {m+1} \left( \E K_{m+1}^- + \frac { m+1-\E K_{m+1}^-} m \E K_m^-
\right) {,}
 $$
 and reorganizing the equation yields the second recursive formula.
 \qed

 \begin{lem}\label{xia4.8}
If $\a_m > \b_m$, then
 $$
1 + \E  K_m^+ \le \frac{\a_m} {\a_m - \b_m }.
 $$
If $\b_m > \a_m$, then,
 $$
\E K_m^- \le \frac  {\b_m}  {\b_m - \a_m}.
 $$
 \end{lem}

\noindent{\it Proof.} Suppose $\a_m > \b_m$. By Lemma~\ref{EK} and that $\E
K_{m-1}^+ \ge 0$, we have
  \begin{equation}\label{xiEKplus}
1 + \E K_m^+  \le 1 +  \frac{ \b_m  }{ \a_m } \big( 1 + \E K_{m-1}^+
\big) , \ \ \ \forall \ m \ge 1.
  \end{equation}
Iterating (\ref{xiEKplus}) and noticing that $\b_k / \a_k$ is
increasing in $k$ as well as $\E K_0^+ = 0$, we conclude that
 $$
1 + \E K_m^+ \le \sum_{i=0}^{m-1} \( \frac{\b_m}{\a_m} \) ^l + \(
\frac{\b_m}{\a_m} \) ^m \big( 1 + \E K_0^+ \big) \le \frac{1} { 1 -
\frac{\b_m}{\a_m} }=  \frac{\a_m} {\a_m - \b_m } .
 $$

Assume $\a_m < \b_m$. Using Lemma~\ref{EK} again together with the fact that $\E K_{m+1}^- \ge 1$, we have
\begin{equation} \label{foriterate}
\E K_m^-  \le 1 + \frac{  \a_m}{ \b_m }  \E K_{m+1}^- .
\end{equation}
Noticing that $\a_k / \b_k$ is decreasing in $k$, we conclude that
 $$
\E K_m^- \le  \sum_{i=0}^{l-1} \( \frac{\a_m}{ \b_m } \) ^i + \(
\frac{\a_m}{ \b_m } \) ^l  \E K_{m+l}^-
 $$
by iterating (\ref{foriterate}). Recalling $\E K_{m+l}^- \le
(m+l)$, we have, by letting $l \rightarrow \infty$, that
 $$
\E K_m^- \le  \sum_{i=0}^\infty \( \frac{\a_m}{ \b_m } \) ^i =
\frac{1}{1 - \frac{\a_m}{\b_m}} = \frac{\b_m}{\b_m - \a_m} .
 $$
  \qed

\section{Proof of Theorem~\ref{ld1}}
\label{prooftheoremld1}
\setcounter{equation}{0}

Let $X$ be a point process with distribution $\bp_{a,b;0;\nu}$, then
by the triangle inequality, we have
$$d_2(\L\Xi,\bp_{a,b;0;\nu})\le d_2(\L\Xi,\L(\m\Xi))+d_2(\L(\m\Xi),\L(\m X))+d_2(\L(\m X),\bp_{a,b;0;\nu}).$$
It follows from \Ref{shuffle1} that both $d_2(\L\Xi,\L(\m\Xi))$ and
$d_2(\L(\m X),\bp_{a,b;0;\nu})$ are bounded by $d_0(\g)$, so it
remains to estimate $d_2(\L(\m\Xi),\L(\m X))$. Clearly, $\m X \sim
\bp_{a,b;0;\nu'}$, where
$$\nu'(dx)=\sum_{i=1}^k\nu(G_i)\d_{t_i}(dx).$$
Using the Stein equation \Ref{Ah} with $\bp=\bp_{a,b;0;\nu'}$, it
suffices to show that for each $f\in\F$,
 \eqa
&&|\E\A h_f(\m\Xi)| = |\E f(\m\Xi)-\bp_{a,b;0;\nu'}(f)| \non
 \\ &\le&
\int_\G\E\left[(1+b)(\epsilon_{1,y}(\Xi_y)+\epsilon_{1,y}(\Xi))+b\bar
r_y(\Xi)\Xi_y(A_y)+b\epsilon_{2,y}(\Xi_y)\right] \l
(dy).\label{ld1-2}
 \ena
To simplify the notation, we fix $f\in\F$, write
$f'(\eta)=f(\m\eta)$, $h'(\eta)=h_f(\m\eta)$ and define $$ \D
h'(\xi;x)=h'(\xi+\d_x)-h'(\xi). $$ Noting that $h'$ acts on the
`shuffled' configurations so one can swop $\nu'$ for $\nu$ in $\A
h'$, we apply \Ref{eq:Palm} to expand $\E\A h'(\Xi)$ as
 \eqa
\E\A h'(\Xi)&=&b\int_\G \int_\G \E [\D h'(\Xi_y+\d_y;x)-\D h'(\Xi;x)]\l(dy)\nu(dx)\non\\
&& + \int_\G \E [-\D h'(\Xi_x;x)+\D h'(\Xi;x)]\l(dx)\non\\
&& + \int_\G \E \D
h'(\Xi;x)[a\nu(dx)+b|\l|\nu(dx)-\l(dx)].\label{ld1-4}
 \ena
The last term vanishes since $(a+b|\l|) \nu = \l$, which is ensured by the facts
that $|\nu| = 1$ and $a=(1-b)|\l|$.

To study the first term in (\ref{ld1-4}), we take a coupling $(\T_y,
\U_y, \Pi_y)$ of $\Xi |_{A_y^c}$ (notice that it has the same
distribution as that of $\Xi_y |_{A_y^c}$), $\Xi |_{A_y}$, and
$\Xi_y |_{A_y}$, such that $\L(\T_y+ \U_y) = \L\Xi$ and $\L(\T_y +
\Pi_y) = \L(\Xi_y)$. Dropping the subscript $y$ from $(\T_y, \U_y,
\Pi_y)$, we can write
 \eqs
& & \E\{\D h'(\Xi_y+\d_y;x)-\D h'(\Xi;x)\}
 \\ & = &
\E\{\D h' (\T + \Pi + \d_y; x) - \D h' (\T + \U; x)\}
 \\ & = &
\E\left\{[\D h' (\T + \Pi + \d_y; x) - \D h' (\T; x)] + [\D h' (\T; x) - \D h'
(\T + \U; x)]\right\}.
 \ens
When expanded telescopically, it is the sum of $|\Pi| + 1$ positive $\D_2 h'$-functions for the term in the first pair of square brackets, and $|\U|$ negative
$\D_2 h'$-functions for the term in the second pair of square brackets. Similarly, the second term in (\ref{ld1-4}) can be expressed as the sum of $|\U|$
positive $\D_2 h'$-functions and $|\Pi|$ negative $\D_2 h'$-functions. Therefore, when
 \eq
b\int_{\G}\E(\Xi_y(A_y)+1-\Xi(A_y))\l(dy)+\int_\G\E
(\Xi(A_y)-\Xi_y(A_y))\l(dy)=0, \label{ld1-8}
 \en
the expected numbers of positive and negative $\D_2 h'$-functions are then balanced. Noting
that
 \eq
\int_{\G}\E(\Xi_y(A_y)-\Xi(A_y))\l(dy)=\var(|\Xi|)-\E|\Xi|,\label{ld1-9}
 \en
we obtain (\ref{ld1-8}) by taking $b=\frac{\var(|\Xi|)-\E|\Xi|}{\var(|\Xi|)}$. Now, we denote $\Pi = \sum_{j=1}^{|\Pi|} \d_{x_j}$, $\U = \sum_{j=1}^{|\U|} \d_{y_j}$, and for $\eta=\sum_{i=1}^n\d_{z_i}$, write
$\langle\eta\rangle_0=0$, $\langle\eta\rangle_j=\sum_{i=1}^j\d_{z_i}$ for $1\le j\le n$. Taking $\hat \Xi$ as an independent copy of $\Xi$, we can expand $\E\A h'(\Xi)$ into
 $$
\E\A h'(\Xi) = e_1 + \dots+ e_5,
 $$
where
 \eqs
 e_1 &=&
b \iint_{\G^2} \E
\sum_{j=1}^{|\Pi|}[\D_2h'(\T+\langle\Pi\rangle_{j-1}+\d_y;x,x_j) -
\E \D_2h'(\hat\Xi;z,z)]\l(dy)\nu(dx),
 \\ e_2 & =&
b \iint_{\G^2} \E [\D_2h'(\T;x,y) - \E
\D_2h'(\hat\Xi;z,z)]\l(dy)\nu(dx),
 \\ e_3 & = &
-b \iint_{\G^2} \E
\sum_{j=1}^{|\U|}[\D_2h'(\T+\langle\U\rangle_{j-1};x,y_j) - \E
\D_2h'(\hat\Xi;z,z)]\l(dy)\nu(dx),
 \\ e_4 & = &
-\int_\G \E
\sum_{j=1}^{|\Pi|}[\D_2h'(\T+\langle\Pi\rangle_{j-1};x,x_j) - \E
\D_2h'(\hat\Xi;z,z)]\l(dx),
 \\ e_5 & = &
\int_\G \E \sum_{j=1}^{|\U|}[\D_2h'(\T+\langle\U\rangle_{j-1};x,y_j)
- \E \D_2h'(\hat\Xi;z,z)]\l(dx).
 \ens

Now we concentrate on estimating $e_1$, since others are similar.
Recalling that $\XiyAyc$ is not independent of $\XiyAy$ while
$\XiyByc$ is, we can extract the part as $\XiyByc$ from $\Theta\sim\L(\XiyAyc)$,
and denote it by $\Theta_1$. Take a more detailed coupling $(\T_1,
\T_2, \U, \Pi)$ such that $(\T_1, \T_2)$ is a coupling of $\XiByc$
and $\Xi |_{B_y \setminus A_y}$ (as well as  $\XiyByc$ and $\Xi_y
|_{B_y \setminus A_y}$), and $\T_1$ is dependent of $(\U, \Pi)$. We
then take $(\hat {\T}_2, \hat \U)$ as a copy of $(\T_2, \U)$ such
that $(\hat {\T}_2, \hat \U)$ is independent of $\Pi$ and
$\L(\T_1+\hat {\T}_2+\hat {\U})=\L\Xi$. We insert $\D_2 h' (\T_1; x,
x_j)$ and $\D_2 h' (\T_1; z,z)$ into the square brackets in $e_1$ to
obtain
 $$
e_1 = b \iint_{\G^2} (e_{11} + \cdots + e_{15}) \l(dy)\nu(dx),
 $$
where
 \eqs
 e_{11} & = &
\E \sum_{j=1}^{|\Pi|} [\D_2h'(\T_1 + \T_2
+\langle\Pi\rangle_{j-1}+\d_y;x,x_j) - \D_2h'(\T_1
+\langle\Pi\rangle_{j-1}+\d_y;x,x_j)], \non
 \\ e_{12} & = &
\E \sum_{j=1}^{|\Pi|} [\D_2h'(\T_1
+\langle\Pi\rangle_{j-1}+\d_y;x,x_j) - \D_2h'(\T_1 +
\d_y;x,x_j)], \non
 \\ e_{13} & = &
\E \sum_{j=1}^{|\Pi|} [\D_2h'(\T_1 + \d_y;x,x_j) - \D_2 h'
(\T_1; x,x_j) ] ,
 \\ e_{14} & = &
\E \sum_{j=1}^{|\Pi|} [ \D_2 h' (\T_1; x, x_j) - \D_2 h' (\T_1; z,z)],
 \\ e_{15} & = &
\E |\Pi| \E [\D_2 h' (\T_1;z,z) - \D_2 h'(\T_1 + \hat
{\T}_2 + \hat \U;z,z)].
 \ens

{\em Estimates of $e_{11}$ and $e_{15}$.} Notice $e_{11}$ can
be further decomposed as
 $$
\E \sum_{j=1}^{|\Pi|} \sum_{i=1}^{|\T_2|} \big[ \D_2 h' (\T_1 + \d_y + \langle \T_2, \Pi \rangle_{i,j-1}; x, x_j) - \D_2 h'
(\T_1 + \d_y + \langle  \T_2, \Pi\rangle_{i-1, j-1}; x, x_j) \big],
 $$
where $\langle  \T_2, \Pi \rangle_{i,j} = \langle \T_2 \rangle_{i}
+ \langle \Pi \rangle_j$ are measurable to $( \T_2, \Pi)$.
When we take the expectation conditional on {$\Xi_y |_{B_y}$},
or equivalently on $(\T_2, \Pi)$, it can be interchanged with the
sums. Therefore, we concentrate on the conditional expectation
 \eq
\E \( \left.  \D_2 h' (\T_1 + \d_y + \langle \T_2,
\Pi\rangle_{i,j-1}; x, x_j) - \D_2 h' (\T_1 + \d_y + \langle
\T_2, \Pi\rangle_{i-1,j-1}; x, x_j) \right| \T_2, \Pi \) .
\label{E.D2hbyr1}
 \en
Since by (\ref{mainresult2}), there is no uniform bound for
$\D_2h'$, we write
$$\D_2h'=h^{(1)}+h^{(2)},$$
where
$$h^{(1)}=\min\left\{\max\left(\D_2h',-\frac{2u+5}a\right),\frac{2u+5}a\right\},\ h^{(2)}=\D_2h'-h^{(1)}.$$
Since $$|\D_2h'(\xi;x,y)|\le \frac{2u+5}a\mbox{ for }1+|\xi|>\frac
au,$$ we have \eq |h^{(1)}|\le\frac{2u+5}{a},\ |h^{(2)}|\le 2,\mbox{
and }h^{(2)}(\xi;x,y)=0\mbox{ for }1+|\xi|>\frac
au.\label{ld1-11}
 \en
For the quantity given in (\ref{E.D2hbyr1}), the differences based
on $h^{(1)}$ and $h^{(2)}$ are respectively bounded by the second
and the first terms of $r_y (\Xi_y )$, recalling that $\Xi_y
|_{B_y}$ is equivalent to {$(\T_2, \Pi)$}. Hence,
 \eq
|e_{11}| \le \E |\Pi|\cdot|\T_2| r_y (\Xi_y) = \E r_y (\Xi_y) \Xi_y
(A_y) \Xi_y (B_y \setminus A_y) . \label{E.e11}
 \en

Similarly, taking conditional expectation on $(\hat \T_2,, \hat \U)$, we
get
 \eq
|e_{15}| \le \E |\Pi| \E \big( |\hat \T_2 + \hat \U| r_y (\T_1+\hat
\T_2 + \hat \U) \big) = \E r_y (\Xi) \Xi (B_y) \E \Xi_y (A_y) .
\label{E.e15}
 \en

{\em Estimates of $e_{12}$ and $e_{13}$.} Notice that $\T_2$ disappears now and $\T_1$ is independent of $\Pi$. We use the
conditional expectation on $\Pi$, and find each conditional
expectation, actually being the mean, is less than $\bar r_y (\Xi_y) = \bar r_y (\Xi)$.
Hence
 \eqa
&&|e_{12}| \le \E \frac {|\Pi| (|\Pi| - 1) }2 \bar r_y(\Xi)= \bar
r_y(\Xi) \E \frac { \Xi_y (A_y) ( \Xi_y (A_y) - 1)}{2},
\label{E.e12}\\
&&|e_{13}| \le \bar r_y (\Xi) \E |\Pi| = \bar r_y(\Xi) \E \Xi_y (A_y).
\label{E.e13}
 \ena

{\em Estimate of $e_{14}$.} In fact, $e_{14}$ is
another kind of difference that is very different from the other four since the two
point processes have the same size. Let us state a result which tells
us the cost of shuffling points $x$ and $y$ in $\D_2h(\xi;x,y)$. Define
 $$
D h' (\xi;x,y) = h'(\xi+\d_x) - h'(\xi+\d_y), \ \ \ D_2
h'(\xi;x,y;z)=Dh'(\xi+\d_z;x,y)-Dh'(\xi;x,y).
 $$
Then, one can directly verify the following equation:
 \eq
\D_2h'(\xi;x,y)-\D_2h'(\xi;z,z)=D_2h'(\xi;y,z;x)+D_2h'(\xi;x,z;z).\label{ld1lemma-1}
 \en
Consequently, we can rewrite
 $$
e_{14} = \E\sum_{j=1}^ {|\Pi|} [D_2h'(\T_1; x_j,z;x)+D_2h'(\T_1;x,z;z)],
 $$
bearing in mind $\Pi = \sum_{j=1}^{|\Pi|} \d_{x_j}$. Now we estimate {$D_2
h'$. Recalling $|Dh'| \le C_n$ defined in (\ref{Cn}) and estimated in (\ref{estimateCn}), we have}
 $$
\vert Dh'(\xi;x,y)\vert\le 1\wedge
\left(\frac1{2(|\xi|+1)}+\frac1a\right).
 $$
If we set
 $$Dh'=h^{(3)}+h^{(4)},$$
where
$$h^{(3)}=\max\left\{\min\left(Dh',\frac{u+2.5}a\right),-\frac{u+2.5}a\right\}\mbox{ and }h^{(4)}=Dh'-h^{(3)},$$
then
 \eq
\vert h^{(3)}\vert\le \frac{u+2.5}a,\ \vert h^{(4)}\vert\le 1\mbox{
and }h^{(4)}(\xi;x,y)=0\mbox{ for }1+|\xi|>\frac au.\label{ld1-16}
 \en
Comparing with (\ref{ld1-11}), we conclude that $\D_2 h'$, as the
difference of $\D h'$, has conditional expectation (that reduces to its
expectation) less than a half of $\bar r_y(\Xi)$. Therefore,
 \eq
|e_{14}| \le \bar r_y(\Xi) \E |\Pi| = \bar r_y(\Xi) \E \Xi_y (A_y).
\label{E.e14}
 \en

Collecting (\ref{E.e11}-\ref{E.e13}) and (\ref{E.e14}), we obtain
 \eqa
|e_1|
 &\le&
b\int_\G[\E r_y(\Xi_y)\Xi_y(A_y)\Xi_y(B_y\setminus A_y)\non\\
&&+\bar r_y(\Xi)\E(\Xi_y(A_y)+3)\Xi_y(A_y)/2+\E
r_y(\Xi)\Xi(B_y)\E\Xi_y(A_y)]\l(dy).\label{ld1-18}
 \ena

The same procedure can be applied to estimate $e_2$ to $e_5$ by
first selecting the `stepping stones' $\XiyByc$ and $\hXiByc$ to
`bridge' $\XiyAyc$ and $\hXi\sim\L\Xi$ for $e_2$ and $e_4$, and
$\XiByc$ and $\hXiByc$ to `bridge' $\XiAyc$ and $\hXi\sim\L\Xi$ in
$e_3$ and $e_5$, then telescoping within the layer of dependence and
using \Ref{ld1lemma-1} and \Ref{ld1-16} to deal with relocation of
points. We omit the details here and the estimates are summarized
below: \eqa
|e_2|&\le& b\int_\G[\E r_y(\Xi_y)\Xi_y(B_y\setminus A_y)+\bar r_y(\Xi)+\E r_y(\Xi)\Xi(B_y)]\l(dy),\non\\
|e_3|&\le&b\int_\G[\E r_y(\Xi)\Xi(A_y)\Xi(B_y\setminus A_y)+\bar r_y(\Xi)\E(\Xi(A_y)+1)\Xi(A_y)/2\non\\
&&\mbox{\hskip1cm}+\E r_y(\Xi)\Xi(B_y)\E\Xi(A_y)]\l(dy),\non\\
|e_4|&\le&\int_\G[\E r_x(\Xi_x)\Xi_x(A_x)\Xi_x(B_x\setminus A_x)+\bar r_x(\Xi)\E(\Xi_x(A_x)+1)\Xi_x(A_x)/2\non\\
&&\mbox{\hskip1cm}+\E r_x(\Xi)\Xi(B_x)\E\Xi_x(A_x)]\l(dx),\non\\
|e_5|&\le&\int_\G[\E r_x(\Xi)\Xi(A_x)\Xi(B_x\setminus A_x)+\bar r_x(\Xi)\E(\Xi(A_x)+1)\Xi(A_x)/2\non\\
&&\mbox{\hskip1cm}+\E r_x(\Xi)\Xi(B_x)\E\Xi(A_x)]\l(dx).\non \ena
Now, the above four estimates, together with (\ref{ld1-18}), yield
\Ref{ld1-2}, completing the proof of Theorem~\ref{ld1}. \qed

\section{Proof of Theorem~\ref{ld2}}
\label{prooftheoremld2}
\setcounter{equation}{0}

The proof is similar to that of Theorem~\ref{ld1} with some modification to suit the estimation involving the second order reduced Palm processes. Let $Y$ be a point process with distribution $\bp_{a,0;\b;\nu}$, it follows from the triangle inequality that
$$d_2(\L\Xi,\bp_{a,0;\b;\nu})\le d_2(\L\Xi,\L(\m\Xi))+d_2(\L(\m\Xi),\L(\m Y))+d_2(\L(\m Y),\bp_{a,0;\b;\nu}).$$
Again, \Ref{shuffle1} implies that $d_2(\L\Xi,\L(\m\Xi))$ and $d_2(\L(\m Y),\bp_{a,0;\b;\nu})$ are bounded by
$d_0(\g)$, so $d_2(\L(\m\Xi),\L(\m Y))$ is the only term to be estimated.

We replace $\bp$ by $\bp_{a,0;\b;\nu'}$ in the Stein equation
\Ref{Ah} with $\nu'(dx)=\sum_{i=1}^k\nu(G_i)\d_{t_i}(dx)$. It is
sufficient to prove
\eqa \E\A h_f(\m\Xi)&\le&\int_\G\E\left(\epsilon_{1,x}(\Xi_x)+\epsilon_{1,x}(\Xi)\right)\l(dx)\non\\
 &&+\b\iint_{\G^2}\E\left(\epsilon_{1,x,y}(\Xi_{xy})+\epsilon_{1,x,y}(\Xi)+\epsilon_{2,x,y}(\Xi_{xy})\right)\lt(dx,dy)\label{ld2-4}\ena
for all $f\in\F$. For the fixed $f\in\F$, we set $f'(\eta)=f(\m\eta)$, $h'(\eta)=h_f(\m\eta)$ and then apply \Ref{eq:Palm} and \Ref{eq:Palm2} to deduce the following expansion
 \eqa
\E\A h'(\Xi)
 &=&
\int_\G \E [-\D h'(\Xi_x;x)+\D h'(\Xi;x)]\l(dx)\non
 \\ &&
+\b \iint_{\G^2} \E [-\D h'(\Xi_{xy}+\d_y;x)+\D
h'(\Xi;x)]\lt(dx,dy)\non
 \\ & &
+ \int_\G\D \E h'(\Xi;x)
\left(a\nu(dx)-\l(dx)-\b\int_{y\in\G}\lt(dx,dy)\right).\label{ld2-5}
 \ena
The last term of \Ref{ld2-5} vanishes because of the definition of $\nu$ in \Ref{ld2-2}, and $\nu(\G)=1$ ensures
that
$$a=|\l|+\b\iint_{\G^2}\lt(dx,dy)=|\l|+\b(\E|\Xi|^2-|\l|).$$
We take $\hXi$ as an independent copy of $\Xi$ which is also independent of all $\Xi_x$'s and $\Xi_{xy}$'s. {Denote the points in $\XiAx$, $\XixAx$, $\XiAxy$, $\XixyAxy$ respectively by $x_j$, $y_j$, $w_j$, $v_j$. Then using} the two types of local dependence, we have
 \eqa
 &&
\E\A h'(\Xi)\non
 \\ &=&
\int_\G \{\E [-\D h'(\Xi_x;x)+\D h'(\XixAxc;x)]+ \E [\D h'(\Xi;x)-\D
h'(\XiAxc;x)]\}\l(dx)\non
 \\ &&
-\b \iint_{\G^2} [\E \D h'(\Xi_{xy};x)- \E \D
h'(\XixyAxyc;x)]\lt(dx,dy)\non
 \\ &&
-\b \iint_{\G^2} [\E \D h'(\Xi_{xy}+\d_y;x)- \E \D
h'(\Xi_{xy};x)]\lt(dx,dy)\non
 \\ &&
+\b \iint_{\G^2} [ \E \D h'(\Xi;x)- \E \D
h'(\XiAxyc;x)]\lt(dx,dy)\non
 \\ &=&
-\int_\G \E \sum_{j=1}^{{|\Xi_x(A_x)|}}
[\D_2h'(\XixAxc+\langle\XixAx\rangle_{j-1};x,{y_j})-\D_2h'(\hXi;x_1,x_1)]\l(dx)\non
 \\ &&
+\int_\G \E \sum_{j=1}^{|\Xi(A_x)|}
[\D_2h'(\XiAxc+\langle\XiAx\rangle_{j-1};x,{x_j})-\D_2h'(\hXi;x_1,x_1)]\l(dx)\non
 \\ &&
-\b \iint_{\G^2} \E \sum_{j=1}^{|\Xi_{xy}(A_{xy})|}
[\D_2h'(\XixyAxyc+\langle\XixyAxy\rangle_{j-1};x,{v_j})-\D_2h'(\hXi;x_1,x_1)]\lt(dx,dy)\non
 \\ &&
-\b\iint_{\G^2}[\D_2h'(\Xi_{xy};x,y)-\D_2h'(\hXi;x_1,x_1)]\lt(dx,dy)\non
 \\ &&
+\b\iint_{\G^2} \E \sum_{j=1}^{|\Xi(A_{xy})|}
[\D_2h'(\XiAxyc+\langle\XiAxy\rangle_{j-1};x,{w_j})-\D_2h'(\hXi;x_1,x_1)]\lt(dx,dy)\non
 \\ &&
-\E \D_2h'(\hXi;x_1,x_1)\left[\int_\G \E (\Xi_x(A_x)-\Xi(A_x))\l(dx)
+ \b\iint_{\G^2} \E (\Xi_{xy}(A_{xy})+1-\Xi(A_{xy}))\lt(dx,dy)
\right]\non
 \\ & =: &
\phi_1+\dots+\phi_6.\label{ld2-6}
 \ena

The term $\phi_6$ becomes 0 if we set
 $$
\int_\G \E (\Xi_x(A_x)-\Xi(A_x))\l(dx)+ \b\iint_{\G^2} \E
(\Xi_{xy}(A_{xy})+1-\Xi(A_{xy}))\lt(dx,dy)=0,
 $$
hence the $\b$ in \Ref{ld2-1} follows from \Ref{ld1-9},
$\iint_{\G^2}\lt(dx,dy)=\E|\Xi|(|\Xi|-1)$ and the following
observation
$$ \iint_{\G^2} \E (\Xi_{xy}(A_{xy})-\Xi(A_{xy}))\lt(dx,dy)= \iint_{\G^2}\E (|\Xi_{xy}|-|\Xi|)\lt(dx,dy)
=\E(|\Xi|-2-|\l|)(|\Xi|-1)|\Xi|.$$

Following the same steps as the estimation of \Ref{ld1-18}, with
`stepping stones' $\XixBxc$ and $\hXiBxc$ for $\phi_1$, $\XiBxc$ and
$\hXiBxc$ for $\phi_2$, $\XixyBxyc$ and $\hXiBxyc$ for $\phi_3$ and
$\phi_4$, and $\XiBxyc$ and $\hXiBxyc$ for $\phi_5$, we obtain \eqa
\phi_1&\le&\int_\G\E\epsilon_{1,x}(\Xi_x)\l(dx);\non\\
\phi_2&\le&\int_\G\E\epsilon_{1,x}(\Xi)\l(dx);\non\\
\phi_3&\le&\b\iint_{\G^2}\E\epsilon_{1,x,y}(\Xi_{xy})\lt(dx,dy);\non\\
\phi_4&\le&\b\iint_{\G^2}\E\epsilon_{2,x,y}(\Xi_{xy})\lt(dx,dy);\non\\
\phi_5&\le&\b\iint_{\G^2}\E\epsilon_{1,x,y}(\Xi)\lt(dx,dy),\non \ena
which, together with \Ref{ld2-6}, in turn imply \Ref{ld2-4}. This
completes the proof of Theorem~\ref{ld2}. \qed

\noindent\textbf{Acknowledgements}

This work was supported by the Belz fund from the University of
Melbourne (AX) and NSFC10901008 and the National Excellent PhD Thesis fund 200722
from Peking University (FZ).


\def\ac{{Academic Press}~}
\def\aap{{Adv. Appl. Prob.}~}
\def\ap{{Ann. Probab.}~}
\def\anap{{Ann. Appl. Probab.}~}
\def\jap{{J. Appl. Probab.}~}
\def\jws{{John Wiley $\&$ Sons}~}
\def\ny{{New York}~}
\def\ptrf{{Probab. Theory Related Fields}~}
\def\sp{{Springer}~}
\def\spa{{Stochastic Processes Appl.}~}
\def\sv{{Springer-Verlag}~}
\def\tpa{{Theory Probab. Appl.}~}
\def\zw{{Z. Wahrsch. Verw. Gebiete}~}



\begin{thebibliography}{99}
\bibitem{Aldous89} {\sc Aldous, D.} (1989). {\em Probability Approximations via the Poisson Clumping Heuristic.\/} Springer, New
York.

\bibitem{AGG89} {\sc Arratia, R., Goldstein, L. \& Gordon, L.} (1989). Two moments suffice for Poisson
approximations: The Chen-Stein method. \emph{\ap}{\bf 17}, 9--25.

\bibitem{Barbour88}
{\sc Barbour, A. D.} (1988). Stein's method and Poisson process
convergence. \emph{\jap}{\bf 25} (A), 175--184.

\bibitem{BB92} {\sc Barbour, A. D. \& Brown, T. C.} (1992). Stein's method and point
process approximation. \emph{\spa} {\bf 43}, 9--31.

\bibitem{BBX98} {\sc Barbour, A. D., Brown, T. C. \& Xia, A.} (1998). Point processes in time and Stein's
method. \emph{Stochastics and Stochastics Reports}~{\bf 65}, 127--151.

\bibitem{BCL} {\sc Barbour, A. D., Chen, L. H. Y. \& Loh, W.} (1992). Compound
Poisson approximation for nonnegative random variables via Stein's
method. \emph{\ap}{\bf 20}, 1843--1866.

\bibitem{BarbourHall84}
{\sc Barbour,  A. D. \& Hall, P.} (1984). On the rate of Poisson
convergence. \emph{Math. Proc. Cambridge Philos. Soc.}~{\bf 95},
473--480.

\bibitem{bhj} {\sc Barbour,  A. D.,  Holst, L. \& Janson, S.} (1992).
{\em Poisson Approximation.\/} Oxford Univ. Press.

\bibitem{BJ}
{\sc Barbour,  A. D. \& Jensen, J. L.} (1989). Local and tail
approximations near the Poisson limit. \emph{Scandinavian Journal of
Statistics}~{\bf 16}, 75--87.

\bibitem{BM} {\sc Barbour, A. D. \& M{\aa}nsson, M.} (2002). Compound Poisson process approximation. \emph{\ap}{\bf 30}, 1492--1537.

\bibitem{BarbourU98}Ê{\sc Barbour, A. D. \& Utev, S.} (1998). Solving the Stein equation in
compound Poisson approximation. \emph{\aap}{\bf 30},  449--475.

\bibitem{BarbourU99} {\sc Barbour, A. D. \& Utev, S.} (1999). Compound Poisson approximation in total
variation. \emph{\spa}{\bf 82},  89--125.

\bibitem{BarbourXia99} \textsc{Barbour, A. D. \& Xia, A.} (1999).
        Poisson Perturbations. \emph{ESAIM: Probab. Stat.}~{\bf 3}, 131--150.

\bibitem{BarbourXia06} \textsc{Barbour, A. D. \& Xia, A.} (2006). Normal approximation for random sums.
\emph{\aap}\textbf{38}, 693--728.

\bibitem{Brown83}  {\sc Brown, T. C.} (1983). Some Poisson approximations using
compensators. \emph{\ap}{\bf
11}, 726--744.

\bibitem{BWX00} {\sc Brown, T. C., Weinberg, G. V. \& Xia, A.} (2000).
Removing logarithms from Poisson process error bounds. \emph{\spa}{\bf 87}, 149--165.

\bibitem{BX01}  {\sc  Brown, T. C. \& Xia, A.} (2001).  Steins method and birth-death
processes. \emph{\ap}{\bf 29}, 1373--1403.

\bibitem{Ceka97} {\sc \Ceka, V.} (1997). Asymptotic expansions in the exponent: A compound Poisson approach. \emph{\aap}{\bf 29},
374--387.

\bibitem{ChenShao04} {\sc Chen, L. H. Y. \& Shao, Q. M.} (2004). Normal
approximation under local dependence. \emph{\ap}{\bf 32}, 1985--2028.

\bibitem{ChenXia} {\sc Chen, L.~H.~Y.  \& Xia, A.} (2004). Stein's method, Palm theory and Poisson process
approximation. \emph{\ap}{\bf 32}, 2545--2569.

\bibitem{EKM97} {\sc Embrechts, P., Kl\"uppelberg, C. \& Mikosch, T.} (1997).
        {\em Modelling extremal events for insurance and finance.\/}
        Springer Verlag, Berlin.

\bibitem{JacodMano} {\sc Jacod, J. \& Mano, P.} (1988). Une evaluation de la distance entre
les lois d'une semimartingale et d'un processus a accroissements independants.
\emph{Stochastics~}{\bf 25}, 87--124.

\bibitem{Kallenberg83} {\sc Kallenberg, O.} (1983).
{\em Random Measures.\/} Academic Press, London.

\bibitem{Kruopis86} {\sc Kruopis, J.} (1986). Precision of approximations of the generalized Binomial distribution by convolutions of Poisson measures. \emph{Lithuanian Math. J.}~{\bf 26}, 37--49.

\bibitem{NikunenValkeila} {\sc Nikunen, M. \& Valkeila, E.} (1991). A Prohorov bound for a Poisson
process  and an arbitrary counting process with some applications.
\emph{Stochastics and  Stochastics Reports~}{\bf 37}, 133--151.

\bibitem{Presman83} {\sc Presman, E. L.} (1983). Approximation of binomial distributions by infinitely divisible ones. \emph{\tpa}{\bf 28}, 393--403.

\bibitem{Rachev} {\sc Rachev, S. T.} (1991). {\em Probability metrics and the Stability of Stochastic Models.\/}
\jws.

\bibitem{Rollin05}
{\sc R\"ollin, A.}~(2005). Approximation of sum of conditionally independent variables by the translated Poisson distribution. \emph{Bernoulli}~\textbf{11}, 1115--1128.

\bibitem{RU04} {\sc Ruzankin, P. S.}~(2004). On the rate of Poisson process approximation to a Bernoulli process. \emph{\jap}{\bf 41}, 271--276.

\bibitem{SchuhmacherXia08} {\sc Schuhmacher, D. \& Xia, A.} (2008). A new metric between distributions of point processes. \emph{\aap}\textbf{40}, 651--672.

\bibitem{WangXia08} {\sc Wang, X. \& Xia, A.} (2008). On negative binomial approximation to $k$-runs. \emph{\jap}{\bf 45}, 456--471.

\bibitem{Xia93} {\sc Xia, A.} (1993). A note on the Prohorov distance between a counting process and
     a Poisson process.
     \emph{Stochastics and Stochastics Reports\/}~\textbf{45}, 61--77.

\bibitem{Xia97} {\sc Xia, A.} (1997). On using the first difference in the Stein-Chen method. \emph{\anap}{\bf 7}, 899--916.

\bibitem{Xia05} {\sc Xia, A.} (2005). Stein's method and Poisson process approximation. In: {\em An Introduction to Stein's Method,\/} Eds. A.~D.~Barbour \& L.~H.~Y.~Chen, World Scientific
Press, Singapore, 115--181.

\bibitem{XiaZhang08} {\sc Xia, A. \& Zhang, F.} (2008). A polynomial birth-death point process approximation to the Bernoulli process. \emph{\spa}{\bf 118}, 1254--1263.

\end{thebibliography}
\end{document}